%% file: bare_jrnl.tex
\documentclass[journal]{IEEEtran}
\ifCLASSINFOpdf
   \usepackage[pdftex]{graphicx}
\else
\fi
\usepackage[font=small]{caption}
\usepackage{amsmath}
\usepackage{amssymb}
\usepackage{amsthm}
\newtheorem{remark}{Remark}
\newtheorem{definition}{Definition}

\usepackage{subcaption} 
\usepackage{booktabs}
\usepackage{algpseudocode}
\usepackage{scalerel}
\usepackage{verbatim}
\usepackage{commath}
\usepackage{stmaryrd}
\usepackage{rotating}

\algdef{SE}{Begin}{End}{\textbf{begin}}{\textbf{end}}

\newcommand{\mat}[1]{\ensuremath{\mathbf{#1}}}                       
\newcommand{\ten}[1]{\ensuremath{\mathbf{\underline{#1}}}}           
\newcommand{\factorisation}[1]{\ensuremath{\llbracket #1 \rrbracket}}

\newcommand{\td}[1]{\factorisation{#1}}
\newcommand{\unfold}[2]{\ensuremath{\,\ten{#1}_{(#2)}}}
\newcommand{\R}{\ensuremath{\mathbb{R}}}
\newcommand{\stack}[2]{\ensuremath{\,\texttt{stack}_{#1}\big(#2\big)}}
\newcommand{\vect}[1]{\ensuremath{\,\texttt{vec}(#1)}}
\newcommand{\frob}[1]{\ensuremath{\left \|  #1 \right \|_F}}

\newcommand{\A}{\ensuremath{\,\mat{A}}}

\newcommand{\tX}{\ensuremath{\,\ten{X}}}

\newcommand{\tY}{\ensuremath{\,\ten{Y}}}

\input{settings-tikz}
\input{settings-algorithm}

\begin{document}

\title{Reducing Computational Complexity of Tensor Contractions via Tensor-Train Networks}

\author{Ilya~Kisil,
        Giuseppe G. Calvi,
        Kriton~Konstantinidis,
        Yao Lei Xu,
        and~Danilo~P.~Mandic,~\IEEEmembership{Fellow,~IEEE}
\thanks{All authors are with the Department of Electrical and Electronic Engineering, Imperial College London, London SW7 2AZ, U.K. (e-mails: \{i.kisil15, giuseppe.calvi15, k.konstantinidis19, yao.xu15, d.mandic\}@imperial.ac.uk) }}

\maketitle


\section{Introduction}


\IEEEPARstart{T}ENSORS (multi-way arrays) and Tensor Decompositions (TDs) have recently received tremendous attention in the data analytics community, due to their ability to mitigate the Curse of Dimensionality associated with modern large-dimensional Big Data \cite{cichocki2016tensor, Mandic2016_2, Mandic2017}. Indeed, through TDs data volume (e.g. the parameter complexity) is reduced from an exponential one to a linear one in the tensor dimensions, which facilitates applications in areas including compression and interpretability of Neural Networks \cite{cohen2016expressive, cichocki2016tensor}, multi-modal learning \cite{cichocki2016tensor}, and completion of Knowledge Graphs \cite{padia2016inferring, nickel2012factorizing}. At the heart of these TD techniques is the Tensor Contraction Product (TCP), which is used to represent even the most unmanageable higher order tensors through a set of small-scale \textit{core tensors} that are inter-connected via TCP operations.

However, a major obstacle towards a more widespread use of tensors is that even relatively compact multilinear operations, such as the Tucker (TKD) and the Tensor Train (TT) decompositions, tend to be difficult to implement and gain insight into; this is due to the overwhelming indexing and cumbersome mathematical formulae. Consequently, tensor decompositions may not immediately offer the same straightforward intuition as their "flat-view" matrix counterparts, as multi-linear operations become progressively harder to visualize with an increasing dimensionality of the underlying problem. These problems will soon become even more acute as the range of "big data" applications increase.

As a remedy, and taking inspiration from quantum physics and quantum chemistry \cite{cichocki2016tensor,orus2019tensor}, diagrammatic (graphical) approaches have been developed to illustrate tensor based mathematical operations.  In this way, tensors and hence also Tensor Networks (TNs), can be represented as a connection of nodes and edges, in a way that resembles mathematical graphs. Such graphical tools provide an intuitive and concise representation of the underlying multilinear operations, while maintaining the underlying mathematical rigour, but are only recently being adopted in data analytics communities. 

\begin{figure}[t]
    \centering
    \begin{subfigure}[b]{\columnwidth}
        \centering
        \scalebox{0.18}{
            \input{./04-tikz/tensor-tn-rep.tex}
        }
        \caption{        
            A graphical approach to data representation. From left to right, top to bottom: order-0 (scalar), order-1 (vector) of dimension $I_1$, order-2 (matrix) of dimension $I_1 \times I_2$ and order-3 tensor of dimension $I_1 \times I_2 \times I_3$.
        }
        \label{fig:tensor_rep}
    \end{subfigure}
    \begin{subfigure}[b]{\columnwidth}
        \centering
        \vspace{10mm}
        \scalebox{0.15}{
            \input{./04-tikz/basic-operations-tn-rep.tex}
        }
        \caption{
            Graphical representation of several fundamental tensor operations in form of a TN. From top to
            bottom: matrix $\mat{X} \in \R^{I_1 \times I_2}$ and its transpose; contraction
            between matrices $\mat{X}\in \R^{I_1 \times I_2}$ and $\mat{Y}\in \R^{J_1 \times
            J_2}$ along the common mode $I_2=J_1$; contraction between tensors
            $\ten{X}\in\R^{I_1 \times I_2 \times I_3}$ and $\ten{B}\in\R^{J_1 \times J_2
            \times J_3}$ along the common mode $I_3 = J_1$; contraction between tensors
            $\ten{X}$ and $\ten{Y}$ along all modes; unfolding of tensor $\ten{X} \in \R^{I_1 \times I_2 \times I_3}$ along its first mode $I_1$.
        }
        \label{fig:basic_operations_rep}
    \end{subfigure}
    \caption{
        Diagrammatic tensor network representation of: (a) Tensors of various orders, and (b) Basic tensor operations.
    }
    \label{fig:tensor_rep2}
\end{figure}
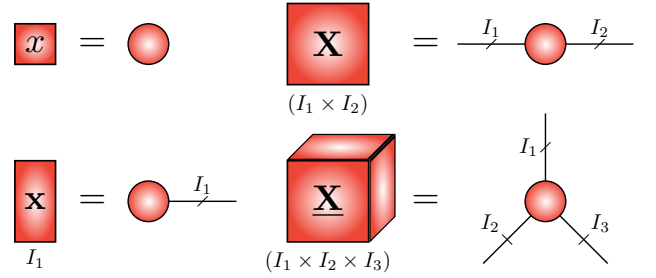
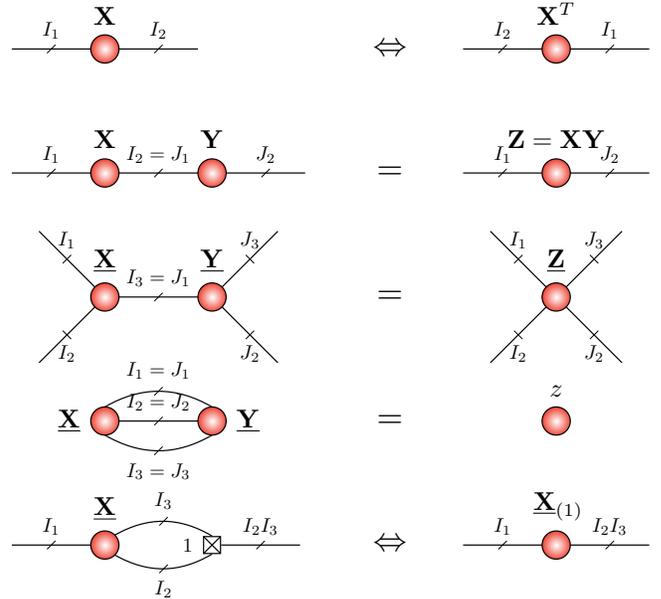

Figure \ref{fig:tensor_rep} shows that a tensor of any order, including scalars (order-0), vectors (order-1), matrices (order-3) and higher order tensors, can be represented by circles (or nodes) in a Tensor Network diagram. The edges of a given tensor represent the tensor modes (dimensions), and, conventionally, they are labelled from the ``top'' of a node (North) in an anti-clockwise direction; the edge labels designate the dimensionality of the mode, while the number of edges connected to the node designates the order of the tensor. A \textit{contraction mode} is an edge that connects two tensor modes of equal size, while an edge which is not connected to another tensor represents a \textit{physical mode} of the underlying tensors, or a TN.

At the core of tensor manipulation is the Tensor Contraction Product (TCP). However, it comes at an exponentially increasing computational complexity of $\mathcal{O}(I^{2N-1})$, with $N$ as the tensor order and $I$ as the modal dimension. To put this into context, for an order-5 tensor with all modal dimensions of $I=10$, the TCP would require a whopping $\mathcal{O}(10^9)$ operations; this is prohibitive to real-world applications based on tensors. The aim of this column is therefore to leverage on the power and flexibility of the graphical tensor network representations in order to develop a Tensor Contraction Product algorithm via Tensor-Train Decomposition (TTD), termed the Tensor-Train Contraction Product (TTCP). This is achieved using purely the diagrammatic representation of tensor networks. Such an approach, in addition to its intuitive derivation, is shown to be computationally efficient, reducing the computational complexity of TCP from the original exponential one in the tensor order to being independent of the tensor order in TTCP. A general and efficient framework for contractions on arbitrary tensor networks on graphs, in the context of spin glass energy and quantum circuits can be found in \cite{pan2020contracting}. It is our hope that this intuitive approach will help better understanding and more efficient implementation of tensor operations, together with promoting their more widespread adoption among researchers and practitioners.



\section{Preliminaries}

\subsection{Basic Tensor Definitions}

\begin{definition}
    \textbf{An order-N tensor}, $\ten{X} \in \R^{I_1 \times \cdots \times I_N}$, is a multi-way array with N modes, with the $n$-th mode of dimensionality $I_n$, for $n=1,\ldots,N$. Special cases of tensors include matrices as order-2 tensors (e.g. $\mat{X} \in \R^{I_1 \times I_2}$), vectors as order-1 tensors (e.g. $\textbf{x} \in \R^{I_1}$), and scalars as order-0 tensors (e.g. ${x} \in \R$).
\end{definition}

\begin{definition}\label{def:mode-n-unfold}
    A \textbf{mode-n unfolding} of a tensor is a procedure of mapping the elements from a
    multidimensional array to a two-dimensional array (matrix). Conventionally, such procedure is associated with stacking mode-$n$ fibers (modal vectors) as column vectors of the resulting matrix.
    For instance, mode-$1$ unfolding of $\ten{X} \in \mathbb{R}^{I_1 \times I_2 \times
    \cdots \times I_N}$ is represented as $\unfold{X}{1} \in \R^{I_1 \times I_2 I_3
    \cdots I_N}$, and given by 

    \begin{equation}
        \unfold{X}{1}\bigg[i_1,\overline{i_2 i_3 \ldots i_N} \bigg] = \tX[i_1,i_2,\ldots, i_N]
    \end{equation}
    Note that the overlined subscripts refer to to linear indexing (or
    Little-Endian) \cite{Dolgov2014}, given by:

    \begin{equation}\label{eq:mode-n-unfold}
    \begin{aligned}
        \overline{i_1 i_2 \dots i_N}
            &= 1 + \sum_{n=1}^N \Bigg[ (i_n - 1) \prod_{n'=1}^{n-1}I_{n'} \Bigg] \\
            &= 1 + i_1 + (i_2 - 1)I_1 + \cdots + (i_n-1)I_1 \ldots I_{N-1}
    \end{aligned}
    \end{equation}
\end{definition}

\begin{definition}
    Any given vector $\mat{x} \in \R^{I_1 I_2 \ldots I_N}$ can be
    \textbf{folded} into an $N$-th order tensor, $\ten{X} \in \R^{I_1 \times I_2 \times \cdots \times I_N}$, with the relation between their entries defined by

    \begin{equation}\label{eq:folding}
        \ten{X}[i_1, i_2, \dots, i_N] = x_i
    \end{equation}
    for all $1\leq i_n \leq I_n$, where $i=1+\sum_{n=1}^{N}(i_n-1)\prod_{k=1}^{n-1}I_k$.
\end{definition}

\begin{definition}
    Consider grouping $J$ order-$N$ tensor samples, $\ten{X} \in \R^{I_1 \times \cdots
    \times I_N}$, so as to form an order-$(N + 1)$ data tensor, $\ten{Y} \in \R^{I_1
    \times \cdots \times I_N \times J}$. This \textbf{stacking} operation is
    denoted by

    \begin{equation}\label{eq:stacking}
        \ten{Y} = \stack{N+1}{\tX[1], \ldots, \tX[J]}
    \end{equation}
    In other words, the so combined tensor samples introduce another dimension, the $(N +
    1)$th mode of $\ten{Y}$, such that the mode-$(N + 1)$ unfolding of $Y$ becomes

    \begin{equation}
        \unfold{Y}{N+1} = \Big[ \vect{\tX[1]}, \ldots, \vect{\tX[J]} \Big]
    \end{equation}
\end{definition}

\begin{definition}
    The \textbf{mode-n product} takes as inputs an order-$N$ tensor, $\tX \in \R^{I_1
    \times I_2 \times \cdots \times I_N}$, and a matrix, $\A \in \R^{J \times I_n}$, to
    produce another tensor, $\tY$, of the same order as the original tensor $\tX$. The
    operation is denoted by
    
    \begin{equation}
        \tY = \tX \times_n \A
    \end{equation}
    where $\tY \in \R^{I_1 \times \cdots \times I_{n-1} \times J \times I_{n+1} \times
    \cdots \times I_N}$. The mode-$n$ product is comprised of 3 consecutive steps:
    
    \begin{equation}
    \begin{aligned}\label{eq:mode-n-product}
        \tX &\rightarrow \unfold{X}{n} \\
        \unfold{Y}{n} &= \A \unfold{X}{n} \\
        \unfold{Y}{n} &\rightarrow \ten{Y}
    \end{aligned}
    \end{equation}
\end{definition}

\begin{definition}
    At the core of TDs and TNs is the \textbf{Tensor Contraction Product} (TCP), an operation similar to mode-$n$ product, but the arguments of which are multidimensional arrays that can be of different order. For instance, given an $N$-th order tensor $\ten{X} \in \mathbb{R}^{I_1\times \cdots \times I_N}$ and another $M$-th order tensor $\ten{Y}\in \mathbb{R}^{J_1\times \cdots \times J_M} $, with common modes $I_n = J_m$, then their $(n,m)$-contraction denoted by $\times^m_n$,  yields a third tensor $\ten{Z}\in \mathbb{R}^{I_1 \times \cdots \times I_{n-1} \times I_{n+1}  \times \cdots \times I_N \times J_1 \times \cdots \times J_{m-1} \times J_{m+1}  \times \cdots \times J_M}$ of order $(N+M-2)$, $\ten{Z}=\ten{X} \times_n^m \ten{Y}$, with entries
    \begin{equation}\label{eq:cont}
    \begin{aligned}
        &z_{i_1,\dots,i_{n-1}, i_{n+1}, \dots, i_N, j_1, \dots, j_{m-1}, j_{m+1}, \dots, j_M   } = \\
            &= \sum_{i_n}^{I_n} x_{i_1, \dots, i_{n-1}, i_n, i_{n+1}, \dots, i_N y_{j_1, \dots, j_{m-1}, i_n, j_{m+1}, \dots, j_M}}
    \end{aligned}
    \end{equation}
\end{definition}

Observe the overwhelming indexing associated with the TCP operation in (8). This quickly becomes unmanageable for larger tensor networks, whereby multiple TCPs are carried out across a large number of tensors. Manipulation of such expressions is prone to errors and prohibitive to manipulation of higher order tensors. This further motivates the use of diagrammatic representations for the execution and understanding of the underlying operations, as illustrated in Figure 2.

\begin{definition}
    The \textbf{Tensor-Train Decomposition} (TTD) is a tensor network that approximates a large order-$N$ tensor, $\ten{X} \in \R^{I_1 \times I_2 \times \cdots \times I_N}$, as a series of $N$ contracting core tensors, $\ten{G}^{(1)}, \dots, \ten{G}^{(N)}$, that is 
    \begin{equation}\label{eq:tt_def}
    \begin{aligned}
        \ten{X} \approx \ten{G}^{(1)} \times_2^1 \ten{G}^{(2)} \times_3^1 \ten{G}^{(3)} \times_3^1 \cdots \times_3^1 \ten{G}^{(N)}
    \end{aligned}
    \end{equation}
\end{definition}

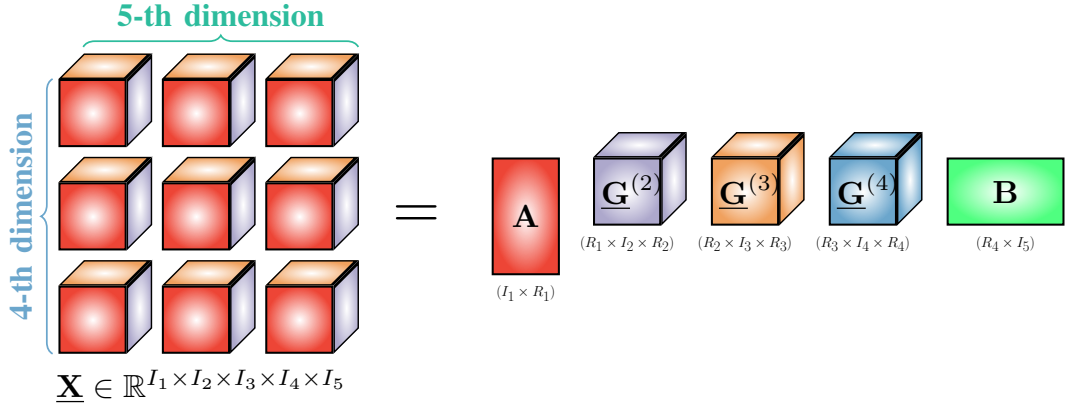
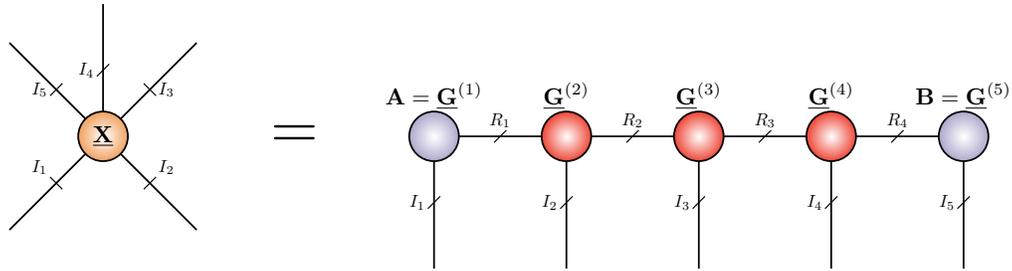
\begin{figure*}[t]
    \centering
    \begin{subfigure}[b]{\textwidth}
        \centering
        \scalebox{0.22}{
            \input{./04-tikz/TT.tex}
        }
        \caption{        
            A Tensor Train type TN decomposition of an order-$5$ tensor, $\ten{X} \in
            \R^{I_1 \times I_2 \times I_3 \times I_4 \times I_5}$, into three order-3 small-scale cores, $\ten{G}^{(2)}, \ten{G}^{(3)}, \ten{G}^{(4)}$, and two order-2 factor matrices, $\mat{A}$, and $\mat{B}$.
        }
        \label{fig:tt1}
    \end{subfigure}
    \begin{subfigure}[b]{\textwidth}
        \centering
        \vspace{10mm}
        \scalebox{0.22}{
            \input{./04-tikz/TT-rep.tex}
        }
        \caption{
            Graphical representation of Tensor Train decomposition of an order-$5$ tensor, $\ten{X} \in \R^{I_1 \times I_2 \times I_3 \times I_4 \times I_5}$. 
        }
        \label{fig:tt2}
    \end{subfigure}
    \caption{
        The Tensor-Train Decomposition (TTD) illustrated in: (a) Conventional way, and (b) Via a tensor network representation.
    }
    \label{fig:tensor_rep3}
\end{figure*}

The cores, $\ten{G}^{(n)}$, within the TTD can be computed through the Tensor-Train Singular Value Decomposition (TT-SVD) algorithm, which is described in Algorithm 1. A TT decomposition of an order-5 tensor is illustrated in Figure \ref{fig:tt1}.
  
\begin{algorithm}
    \caption{Algorithm 1: TT Singular Value Decomposition (TT-SVD) for an N-th order tensor}
    \label{algo:ttsvd}
    
    \Input{
        Data tensor, $\ten{X} \in \R^{I_1 \times I_2 \times \cdots \times I_N}$, and
        desired approximation accuracy, $\epsilon$
    }

    \Output{
        Core tensors, $\ten{G}^{(1)}, \dots, \ten{G}^{(N)}$, approximating $\ten{X} \in
        \R^{I_1 \times I_2 \times \cdots \times I_N}$
    }

    Initialize cores, $\ten{G}^{(1)}, \dots, \ten{G}^{(N)}$, and $R_0=1$

    Compute truncation parameter $\delta = \frac{\epsilon}{\sqrt{N-1}} ||\ten{X}||_F$

    $\ten{Z} \leftarrow \ten{X}$, and $\mat{Z} \leftarrow \mat{Z}_{(1)}$

    \For{$n=1$ to $N-1$}{
        Compute $\delta$-truncated SVD:
        $
            \mat{Z} = \mat{USV}+ \mat{E},
            \text{ s.t. }
            \frob{\mat{E}} \leq \delta; \mat{U} \in \R^{R_{(n-1)}I_n \times R_n}
        $

        $\ten{G}^{(n)} \leftarrow$
            \Reshape{$\mat{U}, [ R_{(n-1)}, I_n, R_n   ]$}

        $\mat{Z} \leftarrow$
            \Reshape{$\mat{SV}^T, [R_n I_{(n+1)}, I_{(n+2)} I_{(n+3)}\dots I_N]) $}
    }

    $\ten{G}^{(N)} \leftarrow \mat{Z}$

    \Return{
         $\ten{G}^{(1)}, \ten{G}^{(2)}, \dots, \ten{G}^{(N)}$
    }    
\end{algorithm}


\subsection{Representing Multilinear Operations via Tensor Networks}

When considered solely through mathematical expressions, the interactions among tensors are extremely exhaustive and lack physical intuition, as is the case with the contraction product in (\ref{eq:cont}). On the other hand, as shown in Figure \ref{fig:basic_operations_rep}, the graphical representation allows for intuitive ways to correctly and intuitively present even the most complex operations.

Moreover, graphical representations may be exploited to dramatically enhance the understanding of the fundamental concepts, owing to a physically meaningful basis of such visualisation, as illustrated in Figure \ref{fig:tt2} on an example of the TTD from Figure \ref{fig:tt1}. This promises to lead to novel approaches to efficient design and implementation of tensor valued algorithms, as shown next on the example of the tensor contraction product.

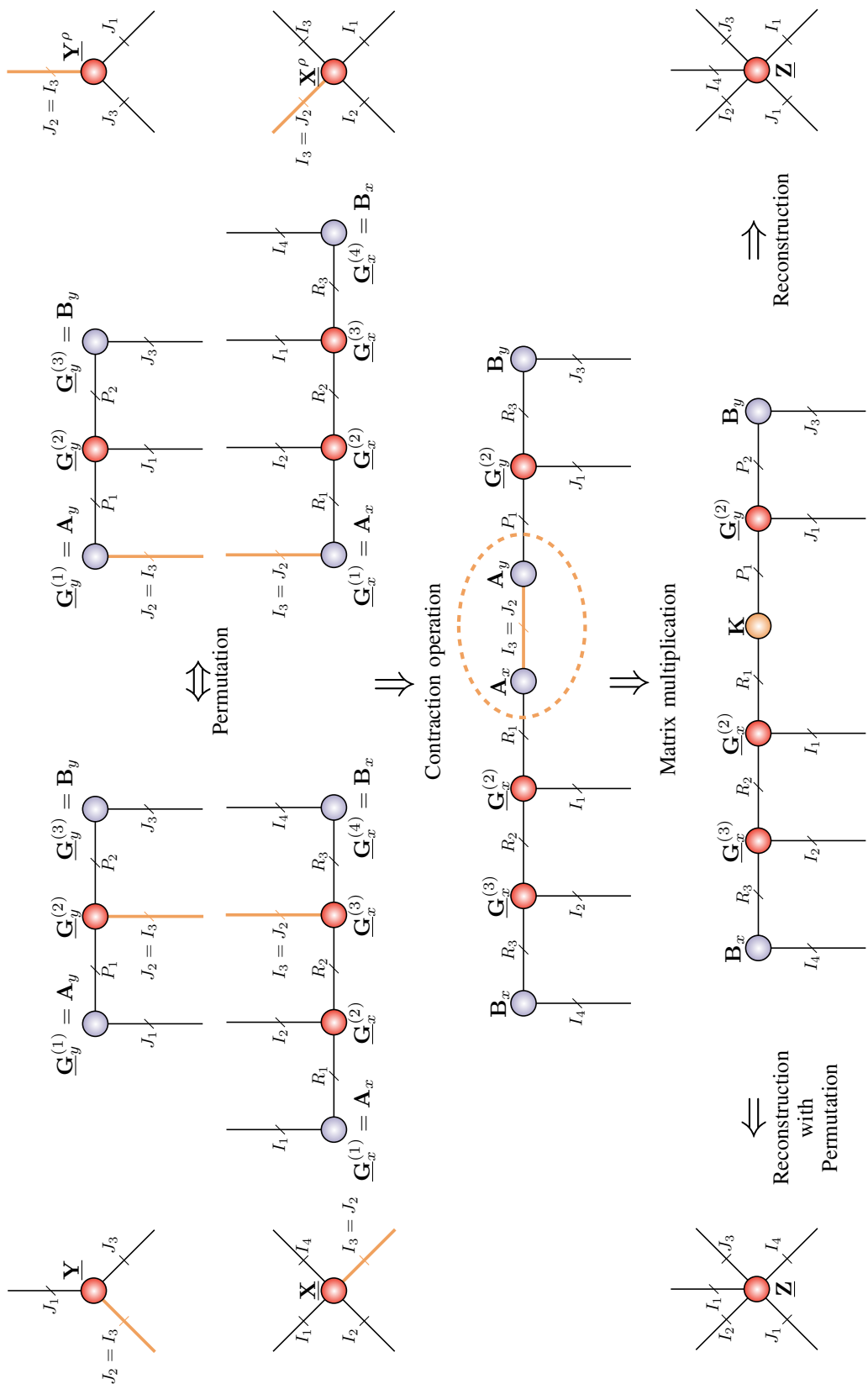
\begin{sidewaysfigure*}
    \centering
    \scalebox{0.2}{
        \input{./04-tikz/TTCP-rep.tex}
    }
    \caption{
        Graphical illustration of the TTCP. The contraction is performed on an order-4 tensor,
        $\ten{X}$, and an order-3 tensor, $\ten{Y}$, along their common modes $I_3 = J_2$. By exploiting the topology of the TT representations of the original tensors, $\ten{X}$ and $\ten{Y}$, the TT decompositions of the permuted versions of $\ten{X}$ and $\ten{Y}$ are obtained, allowing for the reduction of the otherwise computationally prohibitive contraction in (\ref{eq:cont}) to a single matrix multiplication.
    }
    \label{fig:contrpod}
\end{sidewaysfigure*}

\section{Graphical approach to contraction of tensors in the TT-format}\label{sec:ttcp}

We have shown that tensor representation through nodes and edges is both immensely powerful and offers immediate flexibility and versatility in mathematical operations. In this way, it is possible to split links, contract cores, rotate sub-networks, and merge TNs by contracting specific cores when an original tensor is treated as a graph.  In other words, the very topology of a TN is sufficient to obtain insights which would otherwise be difficult to gain and non-obvious.

We next show that the Tensor Contraction Product, when performed graphically through a TT-type of TN, becomes more intuitive, straightforward to implement, and computationally efficient, whereby each core conveys information regarding the corresponding original tensor mode.

\textbf{Example 1.} Given that the first and last cores in a TT are always matrices, the question arises whether it is possible to exploit such an arrangement to reduce an $\times_n^m$ contraction product to a simple matrix multiplication. For illustration, consider a $\times_3^2$ contraction between a 4th-order tensor, $\ten{X}\in \mathbb{R}^{I_1\times I_2 \times I_3 \times I_4}$ and a 3rd-order tensor, $\ten{Y}\in \mathbb{R}^{J_1\times J_2 \times J_3}$, where $I_3=J_2$. This would result in a tensor $\ten{Z}$, given by

\begin{equation}\label{eq:ttcp_cont}
    \ten{X} \times^2_3 \ten{Y} = \ten{Z} \in \mathbb{R}^{I_1 \times I_2 \times I_4 \times J_1 \times J_3}
\end{equation}

The tensors $\ten{X}$ and $\ten{Y}$ can also have their indices permuted to yield
$\ten{X}^\rho \in \mathbb{R}^{I_3 \times I_1 \times I_2 \times I_4}$ and $\ten{Y}^\rho
\in \mathbb{R}^{J_2 \times J_1 \times J_3}$, such that the contracting modes, $I_3$ for
$\ten{X}$, and $J_2$ for $\ten{Y}$ become the first modes. Upon performing the TT decomposition
on $\ten{X}^\rho$ and $\ten{Y}^\rho$ respectively, this yields the constituent cores

\begin{equation}
    \begin{aligned}
        \ten{G}_x^{(n)} &\in \mathbb{R}^{R_{(n-1)} \times I_n \times R_n}, \hspace{2mm} n = 1, 2,3, 4\\
        \ten{G}_y^{(m)} &\in \mathbb{R}^{P_{(m-1)} \times J_m \times P_m}, \hspace{2mm} m=1, 2, 3 \\
        P_0 &= P_3 = 1 \\
        R_0 &= R_4 = 1 \\
    \end{aligned}
\end{equation}

This implies that, since $I_3 = J_2$, the desired contraction in \eqref{eq:ttcp_cont} can
be obtained as a simple matrix product  $\ten{G}_x^{(1)} = \mat{A}_x \in
\mathbb{R}^{I_1 \times R_1}$ and $\ten{G}_y^{(1)} = \mat{A}_y \in \mathbb{R}^{J_2
\times P_1}$, in the form

\begin{equation}
    \mat{A}_x^T\mat{A}_y = \mat{K} \in \mathbb{R}^{R_1 \times P_1}
\end{equation}

The intuition behind Example 1 is at the centre of the proposed Tensor-Train Contraction Product (TTCP), a method to perform contractions between tensors that are represented in the TT-format. This is illustrated in Figure \ref{fig:contrpod}, where the cores $\mat{A}_x$ and $\mat{A}_y$ of the tensors $\ten{X}^\rho$ and $\ten{Y}^\rho$ are linked by a standard matrix multiplication, $ \mat{A}_x^T \mat{A}_y$, connecting the two TTs of $\ten{X}^\rho$ and $\ten{Y}^\rho$. This results in the newly formed TT, whereby the link introduced by the contraction product is consolidated into the matrix, $\mat{K} \in \mathbb{R}^{R_1 \times P_{1}}$. Thus, the task of reconstructing the tensor represented by the union of the two TTs yields $\ten{Z}^\rho \in \mathbb{R}^{I_4 \times I_2 \times I_1 \times J_1 \times J_3}$, which is a permuted version of the tensor $\ten{Z}$ in \eqref{eq:ttcp_cont}. In order to obtain the latter, it is sufficient to permute $\ten{Z}^\rho$ such that the resulting dimensions follow the natural anti-clockwise order obtained by the original contraction. The proposed TTCP procedure is summarised in Algorithm \ref{algo:contprodtt}, where the operator $(\cdot)^{T_{[1;3]}}$ permutes any given tensor such that the first and third modes are exchanged.

\begin{algorithm}
    \caption{Algorithm 2: Tensor-Train Contraction Product (TTCP)}
    \label{algo:contprodtt}
    
    \Input{
        $\ten{X}\in \mathbb{R}^{I_1 \times \cdots \times I_N}$, $\ten{Y} \in \mathbb{R}^{J_1 \times \cdots \times J_M}$, $n$, $m$
    }
    
    \Output{
        $\ten{Z} = \ten{X} \times_n^m \ten{Y}$ 
    }

    \nosemic Perform permutation of an original $\ten{X} \rightarrow \ten{X}^\rho$ such that \;
    \dosemic\nonl $$
        \ten{X}^\rho \in \mathbb{R}^{I_n \times I_1 \times \cdots \times I_{n-1} \times
        I_{n+1} \cdots \times I_N}
    \;$$

    \nosemic Represent $\ten{X}^\rho$ in the TT format using TT-SVD \;
    \dosemic\nonl $$
        \ten{X}^\rho = \mat{A}_x \times_2^1 \ten{G}_x^{(2)} \times _3^1 \cdots \times_3^1
        \ten{G}_x^{(N-1)} \times_3^1 \mat{B}_x
    \;$$

    \nosemic Repeat steps 1 and 2 for $\ten{Y}$ to obtain \;
    \dosemic\nonl $$
        \ten{Y}^\rho = \mat{A}_y \times_2^1 \ten{G}_y^{(2)} \times _3^1 \cdots \times_3^1
        \ten{G}_y^{(M-1)} \times_3^1 \mat{B}_y
    \;$$

    \nosemic Perform contraction $\times_n^m$ through matrix multiplication \;
    \dosemic\nonl $$
        \mat{K} = \mat{A}_x^T \mat{A}_y
    \;$$

    \nosemic Obtain reconstructed version of $\ten{Z}^\rho$ \;
    \dosemic\nonl 
    \begin{equation*}
        \begin{aligned}
        \ten{Z}^\rho =&
            \mat{B}_x \times_2^1 \ten{G}_x^{(N-1)} \times _3^1 
            \cdots \\
            & \cdots \times_3^1 \ten{G}_x^{(2)} \times_3^1 \mat{K} \times_2^1 \ten{G}_y^{(2)} \times _3^1 \cdots \\
            & \cdots \times_3^1 \ten{G}_y^{(M-1)} \times_3^1 \mat{B}_y
        \;        
        \end{aligned}
    \end{equation*}


    \nosemic Perform permutation of $\ten{Z}^\rho \rightarrow \ten{Z}$ such that \;
    \dosemic\nonl $$
\ten{Z} \in \mathbb{R}^{
            I_1 \times \cdots \times I_{n-1} \times I_{n+1} \times \cdots \times I_N \times
            J_1 \times \cdots \times J_{m-1} \times J_{m+1} \times \cdots \times J_M}
    \;$$

    \Return{
        $\ten{Z}$
    }

\end{algorithm}

The benefits of the TTCP are two-fold:

\begin{enumerate}
    \item The resulting tensor $\ten{Z} = \ten{X} \times^2_3 \ten{Y}$ is
        automatically given in the TT-format, and is thus suitable for efficient distributed
        storage, while its cores still contain the physically meaningful modal
        information with respect to the individual tensors $\ten{X}$ and $\ten{Y}$;
    \item It provides a computationally more efficient alternative to the \textit{standard} definition in (\ref{eq:cont}) of the contraction product, especially for higher order tensors.
\end{enumerate}

\section{Computational Complexity of TTCP}

We shall first illustrate the advantages of the graphical TTCP through an example, followed by the analysis of computational complexity.

\begin{figure}
    \centering
    \begin{subfigure}{\columnwidth}
        \centering
        \includegraphics[width=0.9\linewidth]{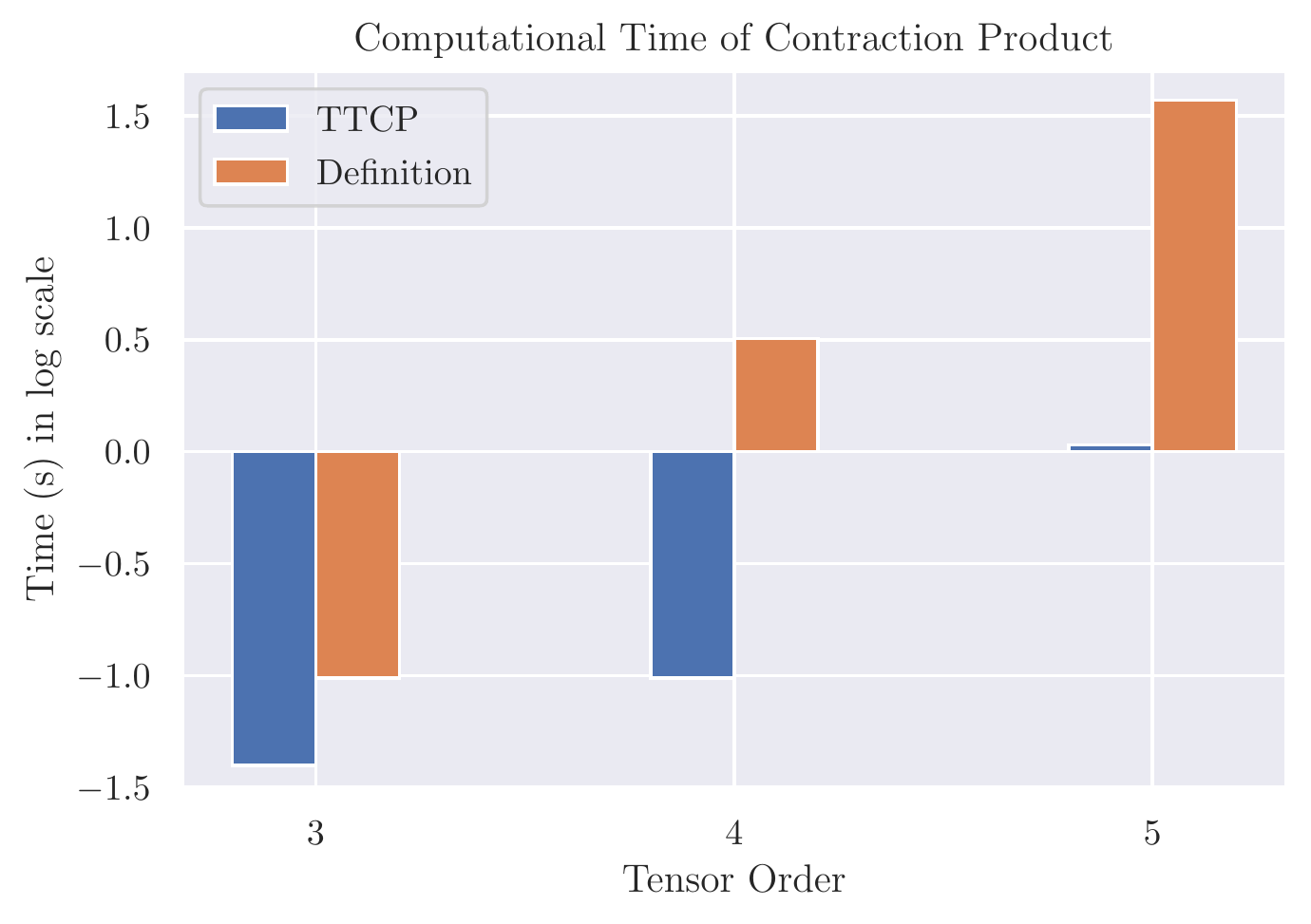}
        \caption{
            Computational time (in seconds, log-scale) between the original expression for tensor contraction in (\ref{eq:cont})
            and the proposed TTCP for the $\ten{X} \times_1^1 \ten{Y}$ contraction. The order-$3$, order-$4$ and
            order-$5$ cases are considered. Stemming from a graphical intuition, the TTCP is
            more computationally efficient than the direct definition of contraction. The computational time of TTCP includes time for both the TTD and the contraction operations.
        }
        \vspace{5mm}
        \label{fig:ttcp}
    \end{subfigure}
    \begin{subfigure}{\columnwidth}
        \centering
        \includegraphics[width=1.0\linewidth]{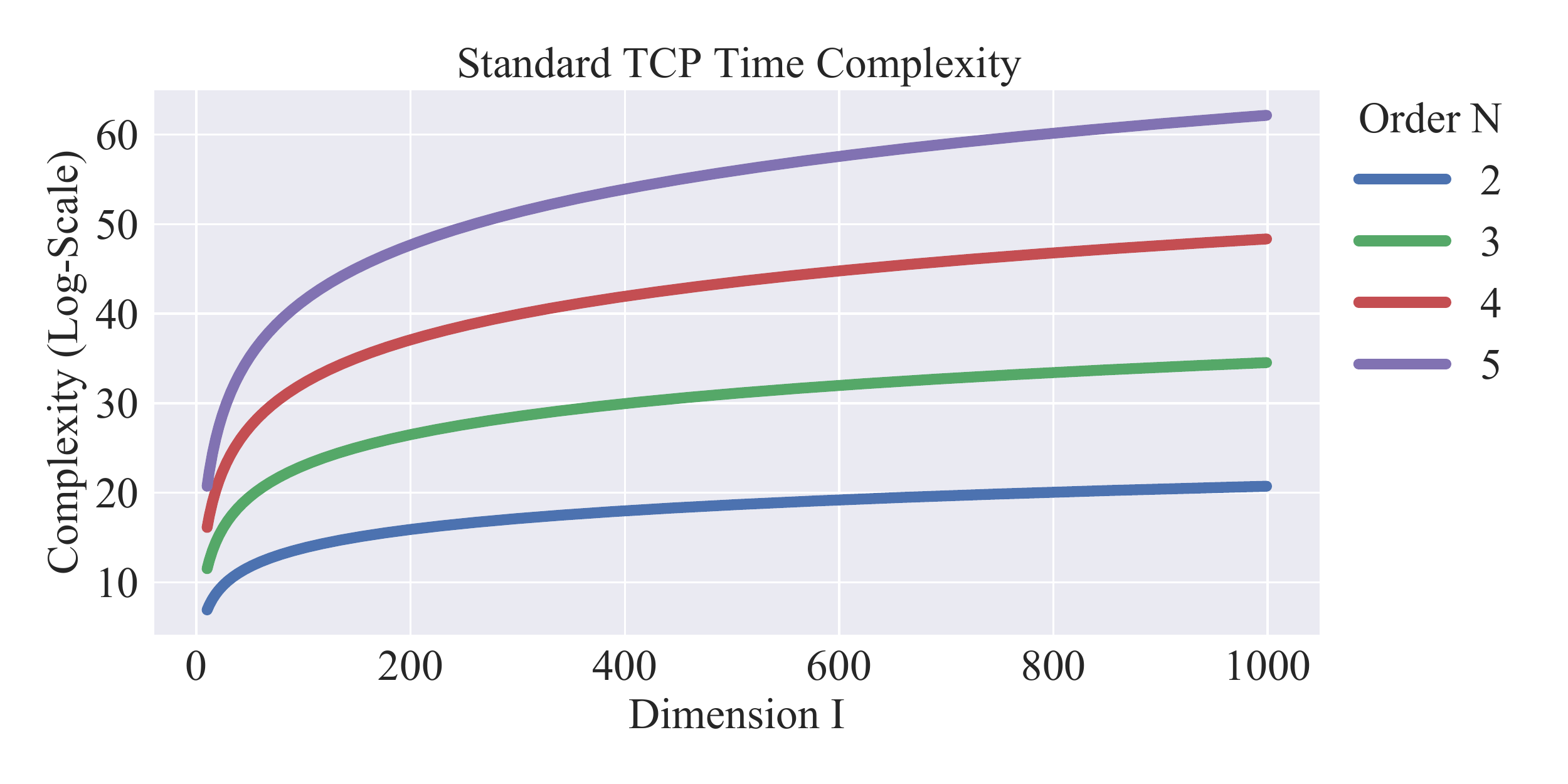}
        \label{fig:tcp_operations1}
    \end{subfigure}
    \begin{subfigure}{\columnwidth}
        \centering
        \includegraphics[width=1.0\linewidth]{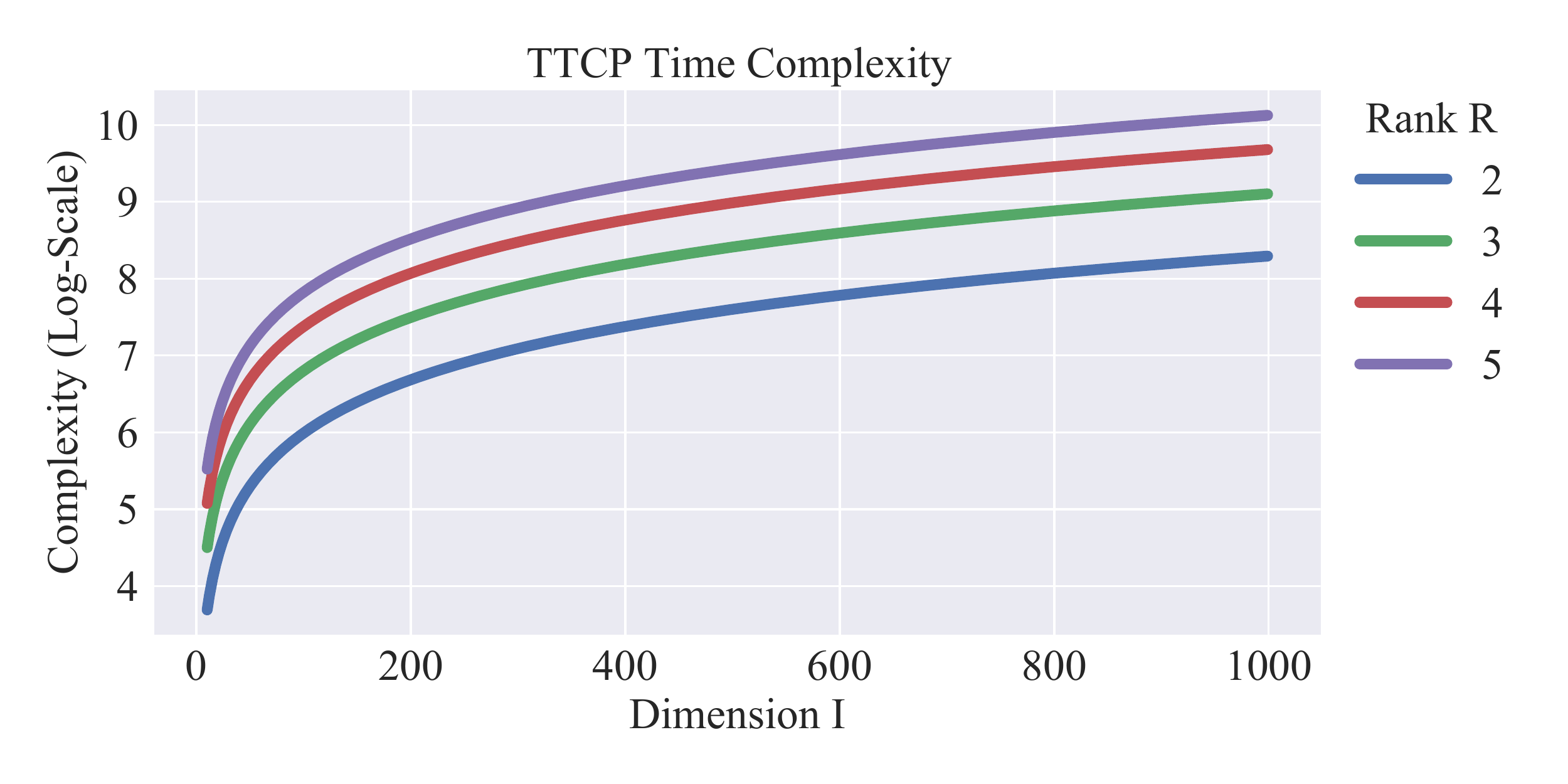}
        \caption{
            Theoretical time complexity (in log-scale), that is, the number of operations required to perform the standard contraction product (top) versus the proposed TTCP (bottom) for varying tensor dimensions, $I$, tensor orders, $N$, and TT ranks, $R$. The theoretical TTCP time complexity, $\mathcal{O}(IR^2)$, is independent of the tensor order and is drastically lower than the complexity of the standard contraction product, $\mathcal{O}(I^{2N-1})$. As noted in Remark 3, if the original tensor is not already in TT format, then it would require an additional operation of TTD, which increases the complexity by $\mathcal{O}(I^{N-1}R^2)$.
        }   
        \label{fig:tcp_operations2}
    \end{subfigure}
    \caption{Complexity analysis of TCP vs TTCP. For illustration, a tensor of order-5 with mode dimensions of $I=1000$ would require $10^{27}$ operations for the standard TCP, while only requiring $25000$ operations for the graphical TTCP with $R=5$, and $4000$ operations for the graphical TTCP with $R=2$.}
    \label{fig:tcp_operations}
\end{figure}

\textbf{Example 2.} The computational efficiency of TTCP is illustrated in Figure \ref{fig:ttcp}, where the two
tensors, $\ten{X}$ and $\ten{Y}$, are drawn from a zero-mean, unit-variance Gaussian distribution, and three different cases are considered: (i) $\ten{X}$, $\ten{Y} \in
\mathbb{R}^{20 \times 20 \times 20}$; (ii) $\ten{X}$, $\ten{Y} \in \mathbb{R}^{20 \times
20 \times 20 \times 5}$; and (iii) $\ten{X}$, $\ten{Y} \in \mathbb{R}^{20 \times 20
\times 20 \times 5 \times 4}$. For each case, the contraction $\ten{X} \times_1^1
\ten{Y}$ was computed both in the original form via \eqref{eq:cont} and also through the proposed TTCP algorithm. The implementations were performed using our own software library HOTTBOX \cite{software-hottbox}, and the computational times were recorded for both cases. Observe that the TTCP exhibits a significant computational efficiency
advantage in all cases, especially for higher order tensors, such as the order-$4$ and order-$5$ tensors considered here, where expression \eqref{eq:cont} tends to become intractable.

\begin{remark}
    The TTCP stems purely from diagrammatic foundations, and makes a full use of the topology
    of the graphical TT representations of tensors $\ten{X}$ and $\ten{Y}$  \cite{10.1137/090752286}.
\end{remark}

\begin{remark}
    The permutation operation in Step 1 of Algorithm 2 needed for TTCP only requires $\mathcal{O}(1)$ operations.
\end{remark}

Generally, instead of performing an expensive $\times_n^m$ tensor contraction product between $\ten{X}\in \mathbb{R}^{I_1 \times \cdots \times I_N}$ and $\ten{Y} \in \mathbb{R}^{J_1 \times \cdots \times J_M}$ according to the standard definition in (\ref{eq:cont}), which incurs $\mathcal{O}(I_1 \cdots I_{n-1} I_{n+1} \cdots I_N I_n J_1 \cdots J_{m-1} J_{m+1} \cdots J_M)$ operations, the proposed TTCP only requires $\mathcal{O}(R_1 I_n P_1)$ operations, where $R_1$ is the first TT rank of $\ten{X}$ and and $P_1$ is the first TT rank of $\ten{Y}$. For illustration, let $I_n=I$ and $J_m=I$ for $n=1, \ldots, N$ and $m=1, \ldots, N$, and $R_1=P_1=R$ in (8). Then, the number of operations required (so called \textit{time complexity}) of the tensor contraction product reduces from the original $\mathcal{O}(I^{2N-1})$, which is \textit{exponential in the tensor order}, $N$, to $\mathcal{O}(IR^2)$ in the TTCP, which is \textit{independent of tensor order}, and is highly efficient for a low TT rank, $R$, as illustrated in Figure \ref{fig:tcp_operations}. 

\begin{remark}
The number of operations required for TTCP is independent of the tensor order, $N$, and is drastically lower compared to the standard definition in (\ref{eq:cont}). The computational advantage of TTCP over standard TCP is $\mathcal{O}(I^{2N-1}) / \mathcal{O}(IR^{2}) \approx \mathcal{O}(I^{2N-3})$. For example, for a mode dimension, $I=1000$, tensor order, $N=5$, and TT rank, $R=5$, the TTCP can reduce the number of operations by almost a factor of $10^{22}$, as shown in Figure \ref{fig:tcp_operations}b in the log-scale. If the original tensor needs first to be decomposed in the TT format, this would require additional computational complexity of $\mathcal{O}(I^{N-1} R^2)$ \cite{9287780}. In this case, the computational advantage of TTCP would still be of the order $\mathcal{O}(I^N)$.
\end{remark}

\section{Conclusion}

There is a significant expansion in both volume and range of applications along with the concomitant increase in the variety of data sources. These ever expanding trends have highlighted the necessity for more versatile analysis tools that offer greater opportunities for algorithmic developments and possibly computationally faster operations than the standard “flat-view” matrix approach.  
Tensors, or multi-way arrays, provide such an algebraic framework which is naturally suited to data of such large volume, diversity, and veracity. Indeed, the associated tensor decompositions have demonstrated, in manifold areas, their potential in breaking the “Curse of Dimensionality” associated with traditional matrix methods, where a necessary exponential increase in data volume leads to adverse or even intractable consequences on computational complexity. A key tool underpinning multi-linear manipulation of tensors and tensor networks is the standard Tensor Contraction Product (TCP). However, depending on the dimensionality of the underlying tensors, the TCP also comes at the price of high computational complexity in tensor manipulation. In this work, we resort to diagrammatic tensor network manipulation to calculate such products in an efficient and computationally tractable manner, by making use of Tensor Train decomposition (TTD). This has rendered the underlying concepts easy to perceive, thereby enhancing intuition of the associated underlying operations, while preserving mathematical rigour. In addition to completely bypassing the cumbersome mathematical multi-linear expressions, and thus eliminating the main source of errors in tensor manipulation, the proposed Tensor Train Contraction Product model is shown to accelerate significantly the underlying computational operations, as it is independent of tensor order and linear in the tensor dimension, as opposed to performing the full computations through the standard approach (exponential in tensor order). It is our hope that this approach will help demystify tensor approaches and encourage researchers and practitioners to explore a whole host of otherwise computationally prohibitive modern applications based on big data and tensors.





\ifCLASSOPTIONcaptionsoff
  \newpage
\fi



\bibliographystyle{IEEEtran}
\bibliography{references}

\end{document}

%% file: settings-tikz.tex
\usepackage[usenames,dvipsnames]{xcolor}
\definecolor{tnode}{rgb}{0.93,0.27,0.21}
\definecolor{mnode}{rgb}{0.68,0.66,0.80}
\definecolor{SeaGreen2}{rgb}{0.31,0.93,0.58}
\definecolor{SkyBlue3}{rgb}{0.42,0.65,0.80}

\definecolor{frontal}{rgb}{0.93,0.27,0.21}
\definecolor{colorFrontal}{rgb}{0.93,0.27,0.21}
\definecolor{lateral}{rgb}{0.68,0.66,0.80}
\definecolor{colorLateral}{rgb}{0.68,0.66,0.80}
\definecolor{horizontal}{rgb}{0.94,0.63,0.37}
\definecolor{colorHorizontal}{rgb}{0.94,0.63,0.37}

\usepackage{tikz}
\usepackage{xifthen}
\usetikzlibrary{
    matrix,
    calc,
    positioning,
    shapes,shapes.geometric,
    arrows,arrows.meta,
    shadows,
    backgrounds,
    decorations.pathreplacing,calligraphy,
    decorations.markings,
    chains,
    mindmap,
    trees
}

\pgfdeclarelayer{background}
\pgfdeclarelayer{foreground}
\pgfsetlayers{background,main,foreground}

\newdimen\XCoord
\newdimen\YCoord
\newdimen\dummyCoord
\newlength\nodelength
\newlength\nodeheight
\newlength\nodeXshift
\newlength\figureHeight

\newcommand*{\ExtractCoordinateX}[1]{
    \path (#1);
    \pgfgetlastxy{\XCoord}{\dummyCoord};
}
\newcommand*{\ExtractCoordinateY}[1]{
    \path (#1);
    \pgfgetlastxy{\dummyCoord}{\YCoord};
}

\newcommand*{\setCoordinateX}[2]{
    \path (#1); \pgfgetlastxy{#2}{\dummyCoord};
}
\newcommand*{\setCoordinateY}[2]{
    \path (#1); \pgfgetlastxy{\dummyCoord}{#2};
}

\makeatletter
    \newcommand\getwidthofnode[2]{
        \pgfextractx{#1}{\pgfpointanchor{#2}{east}}
        \pgfextractx{\pgf@xa}{\pgfpointanchor{#2}{west}}
        \addtolength{#1}{-\pgf@xa}
    }

    \newcommand\getheightofnode[2]{
        \pgfextracty{#1}{\pgfpointanchor{#2}{south}}
        \pgfextracty{\pgf@xa}{\pgfpointanchor{#2}{north}}
        \addtolength{#1}{-\pgf@xa}
    }
\makeatother

\tikzset{
    invisible/.style = {
        opacity=0
    },
    visible on/.style = {alt={#1{}{invisible}}},
    alt/.code args = {<#1>#2#3}{ \alt<#1>{\pgfkeysalso{#2}}{\pgfkeysalso{#3}}},
    arrow/.style = {-{Stealth[scale=1.5]}, rounded corners, very thick},
    circle dotted/.style={
        dash pattern=on 0.5pt off 10pt,
        line cap=round,
        line width=2pt
    },
    mydash/.style={
        line width=3pt,
        dashed,
        color=SeaGreen2,
        dash pattern=on 15pt off 15pt,
    },
    every node/.style = {
        anchor=south west
    },
    my_matrix/.style = {
        matrix,
        fill=white,
        anchor=west,
        left delimiter={[},
        right delimiter={]},
        row sep=-1pt,
        column sep=4mm,
        outer sep=0pt
    },
    my_matrix2/.style = {
        matrix,
        fill=white,
        anchor=west,
        row sep=-1pt,
        column sep=4mm,
        outer sep=0pt
    },
    no sep/.style={
        inner sep=0pt,
        outer sep=0pt
    },
    double line/.style={
        line width=5pt,
        double distance=10pt,
    },
    cross/.style={
        path picture={
            \draw[black]
                (path picture bounding box.south east) -- (path picture bounding box.north west)
                (path picture bounding box.south west) -- (path picture bounding box.north east);
        }
     },
    kronecker/.style={
        draw,
        line width=3pt,
        fill=white,
        circle,
        cross,
        minimum size=1cm,
        anchor=center
    },
    unfolding/.style={
        line width=3pt,
        draw,
        fill=white,
        cross,
        minimum size=0.8cm,
        anchor=center
    },
    round corners/.style={
        line width=3pt,
        rounded corners=5mm
    },
}

\pgfkeys{
    /tensorParts/.is family, /tensorParts,
    default/.style={
        width =3,
        height=3,
        depth=3,
        start at={0,0},
        xshift=0,
        yshift=0,
        zshift=0,
        color frontal = red,
        color lateral = blue,
        color horizontal = yellow,
        color = red,
        inner color = white,
        anchor = south west,
        name = dummy,
        line width = 0.01mm
    },
    width/.estore in = \mywidth,
    height/.estore in = \myheight,
    depth/.estore in = \mydepth,
    start at/.estore in = \mystartat,
    xshift/.estore in = \myxshift,
    yshift/.estore in = \myyshift,
    zshift/.estore in = \myzshift,
    color frontal/.estore in = \mycolorfrontal,         
    color lateral/.estore in = \mycolorlateral,         
    color horizontal/.estore in = \mycolorhorizontal,   
    color/.estore in = \mycolor,                        
    inner color/.estore in = \myinnercolor,
    anchor/.estore in = \myanchor,                      
    name/.estore in = \myname,
    line width/.estore in = \mylinewidth,
    /Tensor/.is family, /Tensor,
    defaultTensor/.style={
        width=3,
        height=3,
        depth=3,
        start at={0,0},
        xshift=0,
        yshift=0,
        zshift=0,
        colorf=colorFrontal,
        colorl=colorLateral,
        colorh=colorHorizontal,
        color=white,
        inner color=white,
        fill opacity=1,
        anchor=south west,
        name=dummy,
        line width=3.5pt,
        on slide=1-
    },
    width/.estore in = \width,
    height/.estore in = \height,
    depth/.estore in = \depth,
    start at/.estore in = \startat,
    xshift/.estore in = \xshift,
    yshift/.estore in = \yshift,
    zshift/.estore in = \zshift,
    colorf/.estore in = \tcolorf,             
    colorl/.estore in = \tcolorl,             
    colorh/.estore in = \tcolorh,             
    color/.estore in = \tcolor,               
    inner color/.estore in = \innercolor,
    anchor/.estore in = \anchor,             
    fill opacity/.estore in = \fillopacity,
    name/.estore in = \name,
    line width/.estore in = \linewidth,
    on slide/.estore in = \onslide,
    /SimplePartsTNs/.is family, /SimplePartsTNs,
    defaultSimplePartsTNs/.style={
        inner color=white,
        tnode color=tnode,
        tnode size=3cm,
        tnode anchor=center,
        mnode color=mnode,
        mnode size=2cm,
        mnode anchor=center,
        leaf size=5cm,
        leaf width=3pt,
        bar size=30pt,
        bar width=2pt,
        raise label=15pt
    },
    inner color/.estore in = \innercolor,
    tnode color/.estore in = \tnodecolor,
    tnode size/.estore in = \tnodesize,
    tnode anchor/.estore in = \tnodeanchor,
    mnode color/.estore in = \mnodecolor,
    mnode size/.estore in = \mnodesize,
    mnode anchor/.estore in = \mnodeanchor,
    leaf size/.estore in = \leafsize,
    leaf width/.estore in = \leafwidth,
    bar size/.estore in = \barsize,
    bar width/.estore in = \barwidth,
    raise label/.estore in = \raiselabel,
    /TNs/.is family, /TNs,
    defaultTNs/.style={
        name=dummy,
        start at={0,0},
        xshift=0,
        yshift=0,
        label={},
        type=e,
        layer=foreground,
    },
    name/.estore in = \name,
    start at/.estore in = \startat,
    xshift/.estore in = \xshift,
    yshift/.estore in = \yshift,
    label/.estore in = \label,
    type/.estore in = \type,
    layer/.estore in = \layer,
}

\newcommand\Tensor[2][]{
    \pgfkeys{/Tensor, defaultTensor, #1}

    \node (\name) at ($(\startat)+(\xshift,\yshift,\zshift)$) [rectangle, draw, line width=\linewidth, minimum width=\width cm, minimum height=\height cm, inner color=\innercolor, outer color=\tcolorf, anchor=\anchor, #2] {};

    \begin{scope}
        \clip ($(\name.north west)+(0,0,0)$) -- ($(\name.north west)+(0,0,-\depth)$) --
        ($(\name.north west)+(\width+\linewidth+\linewidth,0,-\depth)$) -- ($(\name.north
        west)+(\width+\linewidth+\linewidth,0,0)$) -- cycle;
        \draw[inner color=\innercolor, outer color=\tcolorh, #2, line width=1.4*\linewidth] ($(\name.north west)+(0,0,0)$) -- ($(\name.north west)+(0,0,-\depth)$) -- ($(\name.north west)+(\width,0,-\depth)$) -- ($(\name.north west)+(\width,0,0)$) -- cycle;
    \end{scope}

    \begin{scope}
        \clip ($(\name.south east)+(0,0,0)$) -- ($(\name.south east)+(0,0,-\depth)$) --
        ($(\name.south east)+(0,\height+\linewidth,-\depth)$) -- ($(\name.south
        east)+(0,\height+\linewidth,0)$) -- cycle;
        \draw[inner color=\innercolor, outer color=\tcolorl, #2, line width=1.4*\linewidth] ($(\name.south east)+(0,0,0)$) -- ($(\name.south east)+(0,0,-\depth)$) -- ($(\name.south east)+(0,\height,-\depth)$) -- ($(\name.south east)+(0,\height,0)$) -- cycle;
    \end{scope}
}

\newcommand\fslice[2][]{
    \pgfkeys{/Tensor, defaultTensor, #1}
    \node (\name) at ($(\startat)+(\xshift,\yshift,\zshift)$) [rectangle, draw, line width=\linewidth, minimum width=\width cm, minimum height =\height cm, inner color=\innercolor, outer color=\tcolor, anchor=\anchor, #2] {};
}

\newcommand\hslice[2][]{
    \pgfkeys{/Tensor, defaultTensor, #1}
    \begin{scope}
        \coordinate (dummy) at ($(\startat)+(\xshift,\yshift,\zshift)$);
        \clip (dummy) -- ($(dummy)+(0,0,-\depth)$) -- ($(dummy)+(\width,0,-\depth)$) -- ($(dummy)+(\width,0,0)$) -- cycle;
        \draw[inner color=\innercolor, outer color=\tcolor, #2, fill opacity=\fillopacity, line width=\linewidth] (dummy) -- ($(dummy)+(0,0,-\depth)$) -- ($(dummy)+(\width,0,-\depth)$) -- ($(dummy)+(\width,0,0)$) -- cycle;
    \end{scope}
    \coordinate (\name) at ($(dummy)+(\width,0,0)$); 
}

\newcommand\setTNsTikzset[1][]{
    \pgfkeys{/SimplePartsTNs, defaultSimplePartsTNs, #1}
    \tikzset{
        tnode/.style={
            circle,
            no sep,
            inner color=\innercolor,
            outer color=\tnodecolor,
            minimum size=\tnodesize,
            anchor=\tnodeanchor,
            draw=black,
            line width=3pt,
        },
        tnode-cross-line/.style={
            rectangle,
            no sep,
            minimum height=0pt,
            minimum width=\tnodesize,
            draw,
            line width=2pt,
            rotate=-45,
            anchor=center
        },
        mnode/.style={
            circle,
            no sep,
            inner color=\innercolor,
            outer color=\mnodecolor,
            minimum size=\mnodesize,
            anchor=\mnodeanchor,
            draw=black,
            line width=3pt,
        },
        leaf-c/.style={
            rectangle,
            draw,
            no sep,
            anchor=center,
            line width=\leafwidth,
            minimum width=\leafsize,
            minimum height=0pt,
        },
        leaf-s/.style={
            rectangle,
            draw,
            no sep,
            anchor=north,
            line width=\leafwidth,
            minimum width=0pt,
            minimum height=\leafsize
        },
        leaf-n/.style={
            rectangle,
            draw,
            no sep,
            anchor=south,
            line width=\leafwidth,
            minimum width=0pt,
            minimum height=\leafsize
        },
        leaf-e/.style={
            rectangle,
            draw,
            no sep,
            anchor=west,
            line width=\leafwidth,
            minimum width=\leafsize,
            minimum height=0pt
        },
        leaf-w/.style={
            rectangle,
            draw,
            no sep,
            anchor=east,
            line width=\leafwidth,
            minimum width=\leafsize,
            minimum height=0pt
        },
        leaf-se/.style={
            rectangle,
            draw,
            no sep,
            anchor=west,
            line width=\leafwidth,
            minimum width=\leafsize,
            minimum height=0pt,
            rotate=-45
        },
        leaf-sw/.style={
            rectangle,
            draw,
            no sep,
            anchor=east,
            line width=\leafwidth,
            minimum width=\leafsize,
            minimum height=0pt,
            rotate=45
        },
        leaf-ne/.style={
            rectangle,
            draw,
            no sep,
            anchor=south,
            line width=\leafwidth,
            minimum width=0pt,
            minimum height=\leafsize,
            rotate=45
        },
        leaf-nw/.style={
            rectangle,
            draw,
            no sep,
            anchor=south,
            line width=\leafwidth,
            minimum width=0pt,
            minimum height=\leafsize,
            rotate=-45
        },
        dim-c/.style={
            anchor=south,
            inner sep=\raiselabel
        },
        dim-w/.style={
            anchor=south,
            inner sep=\raiselabel
        },
        dim-s/.style={
            anchor=east,
            inner sep=\raiselabel
        },
        dim-n/.style={
            anchor=east,
            inner sep=\raiselabel
        },
        dim-e/.style={
            anchor=south,
            inner sep=\raiselabel
        },
        dim-se/.style={
            anchor=south west,
            inner sep=\raiselabel
        },
        dim-sw/.style={
            anchor=south east,
            inner sep=\raiselabel
        },
        dim-ne/.style={
            anchor=east,
            inner sep=\raiselabel
        },
        dim-nw/.style={
            anchor=west,
            inner sep=\raiselabel
        },
        bar-c/.style={
            rectangle,
            draw,
            no sep,
            anchor=center,
            minimum height=0pt,
            rotate=45,
            minimum width=\barsize,
            line width=\barwidth
        },
        bar-s/.style={
            rectangle,
            draw,
            no sep,
            anchor=center,
            minimum height=0pt,
            rotate=45,
            minimum width=\barsize,
            line width=\barwidth
        },
        bar-n/.style={
            rectangle,
            draw,
            no sep,
            anchor=center,
            minimum height=0pt,
            rotate=45,
            minimum width=\barsize,
            line width=\barwidth
        },
        bar-e/.style={
            rectangle,
            draw,
            no sep,
            anchor=center,
            minimum height=0pt,
            rotate=45,
            minimum width=\barsize,
            line width=\barwidth
        },
        bar-w/.style={
            rectangle,
            draw,
            no sep,
            anchor=center,
            minimum height=0pt,
            rotate=45,
            minimum width=\barsize,
            line width=\barwidth
        },
        bar-nw/.style={
            rectangle,
            draw,
            no sep,
            anchor=center,
            minimum height=0pt,
            rotate=45,
            minimum width=\barwidth,
            line width=\barsize
        },
        bar-ne/.style={
            rectangle,
            draw,
            no sep,
            anchor=center,
            minimum height=0pt,
            rotate=-45,
            minimum width=\barwidth,
            line width=\barsize
        },
        bar-sw/.style={
            rectangle,
            draw,
            no sep,
            anchor=center,
            minimum height=0pt,
            rotate=-45,
            minimum width=\barsize,
            line width=\barwidth
        },
        bar-se/.style={
            rectangle,
            draw,
            no sep,
            anchor=center,
            minimum height=0pt,
            rotate=45,
            minimum width=\barsize,
            line width=\barwidth
        },
    }
}
\newcommand\tnode[3][]{
    \pgfkeys{/TNs, defaultTNs, #1}
    \node[tnode, #2] (\name) at ($(\startat)+(\xshift,\yshift)$){#3}
}

\newcommand\mnode[3][]{
    \pgfkeys{/TNs, defaultTNs, #1}
    \begin{pgfonlayer}{\layer}
        \node[mnode, #2] (\name) at ($(\startat)+(\xshift,\yshift)$){#3};
    \end{pgfonlayer}
}

\newcommand\tnLeaf[4][]{
    \pgfkeys{/TNs, defaultTNs, #1}
        \begin{pgfonlayer}{background}
        \node[leaf-\type, #2] (\name) at ($(\startat)+(\xshift,\yshift)$){};
        \node[bar-\type, #3] () at (\name){};
        \node[dim-\type, #4] () at (\name){\label};
        \end{pgfonlayer}
}

\newcommand\tnLeafCC[5][]{
    \ifthenelse{\isempty{#5}}%
    {\tnLeaf[type=c, #1]{#2}{#3}{#4}}%
    {\tnLeaf[type=c, #1]{minimum width=#5, #2}{#3}{#4}}%
}

\newcommand\tnLeafEE[5][]{
    \ifthenelse{\isempty{#5}}%
    {\tnLeaf[type=e, #1]{#2}{#3}{#4}}%
    {\tnLeaf[type=e, #1]{minimum width=#5, #2}{#3}{#4}}%
}

\newcommand\tnLeafWW[5][]{
    \ifthenelse{\isempty{#5}}%
    {\tnLeaf[type=w, #1]{#2}{#3}{#4}}%
    {\tnLeaf[type=w, #1]{minimum width=#5, #2}{#3}{#4}}%
}

\newcommand\tnLeafNN[5][]{
    \ifthenelse{\isempty{#5}}%
    {\tnLeaf[type=n, #1]{#2}{#3}{#4}}%
    {\tnLeaf[type=n, #1]{minimum height=#5, #2}{#3}{#4}}%
}

\newcommand\tnLeafNE[5][]{
    \ifthenelse{\isempty{#5}}%
    {\tnLeaf[type=ne, #1]{#2}{#3}{#4}}%
    {\tnLeaf[type=ne, #1]{minimum height=#5, #2}{#3}{#4}}%
}

\newcommand\tnLeafNW[5][]{
    \ifthenelse{\isempty{#5}}%
    {\tnLeaf[type=nw, #1]{#2}{#3}{#4}}%
    {\tnLeaf[type=nw, #1]{minimum height=#5, #2}{#3}{#4}}%
}

\newcommand\tnLeafSS[5][]{
    \ifthenelse{\isempty{#5}}%
    {\tnLeaf[type=s, #1]{#2}{#3}{#4}}%
    {\tnLeaf[type=s, #1]{minimum height=#5, #2}{#3}{#4}}%
}

\newcommand\tnLeafSE[5][]{
    \ifthenelse{\isempty{#5}}%
    {\tnLeaf[type=se, #1]{#2}{#3}{#4}}%
    {\tnLeaf[type=se, #1]{minimum width=#5, #2}{#3}{#4}}%
}

\newcommand\tnLeafSW[5][]{
    \ifthenelse{\isempty{#5}}%
    {\tnLeaf[type=sw, #1]{#2}{#3}{#4}}%
    {\tnLeaf[type=sw, #1]{minimum width=#5, #2}{#3}{#4}}%
}

%% file: settings-algorithm.tex
\usepackage[longend, ruled, vlined]{algorithm2e}
\usepackage{cleveref}
\NoCaptionOfAlgo
\LinesNumbered
\SetKwInOut{Input}{Input}
\SetKwInOut{Output}{Output}
\SetKwFunction{NonNegLeastSq}{NonNegLeastSq}
\SetKwFunction{Reshape}{Reshape}

\makeatletter
    \newcommand{\nosemic}{\renewcommand{\@endalgocfline}{\relax}}
    \newcommand{\dosemic}{\renewcommand{\@endalgocfline}{\algocf@endline}}
    \let\oldnl\nl
    \newcommand{\nonl}{\renewcommand{\nl}{\let\nl\oldnl}}
\makeatother

\Crefname{algocf}{Algorithm}{Algorithms}

%% file: 04-tikz/tensor-tn-rep.tex

\newcommand{\one}{3}
\newcommand{\I}{6}
\newcommand{\J}{6}
\newcommand{\K}{5}
\newcommand{\fontOne}[1]{\scalebox{7.5}{#1}}
\newcommand{\fontTwo}[1]{\scalebox{4.5}{#1}}

\begin{tikzpicture}
    \setTNsTikzset

    \fslice[start at={0,0},
            width=\one,
            height=\one,
            anchor=center,
            color=frontal,
            name=scalar]{};
    \node[anchor=west] (equal1) at ($(scalar.east) + (1.5,0)$) {\fontOne{$=$}};
    \node[tnode, outer color=frontal, anchor=west] (nodescalar) at ($(equal1.east) + (1.5,0)$){};
    \node[anchor=center] (xscalar) at  ($(scalar.center) + (0,0)$) {\fontOne{$x$}};

    \fslice[start at=scalar.south,
            yshift=-7,
            name=a1,
            width=\one,
            height=\I,
            anchor=north,
            color=frontal,
            name=vector]{};
    \node[anchor=west] (equal2) at ($(vector.east) + (1.5,0)$) {\fontOne{$=$}};
    \node[tnode, outer color=frontal, anchor=west] (nodevector) at ($(equal2.east) + (1.5,0)$){};
    \node[anchor=center] (xvector) at  ($(vector.center) + (0,0)$) {\fontOne{$\mat{x}$}};
    \tnLeafEE[start at=nodevector.east, label={\fontTwo{$I_1$}}]{}{}{}{};
    \node[anchor=north] (dimsvec) at  ($(vector.south) + (0,-0.4)$) {\fontTwo{$I_1$}};

    \fslice[start at=scalar.east,
            xshift=17,
            width=\J,
            height=\I,
            anchor=west,
            color=frontal,
            name=matrix]{};
    \node[anchor=west] (equal3) at ($(matrix.east) + (3,0)$) {\fontOne{$=$}};
    \node[tnode, outer color=frontal, anchor=west] (nodematrix) at ($(equal3.east) + (6.2,0)$){};
    \node[anchor=center] (xmatrix) at  ($(matrix.center) + (0,0)$) {\fontOne{$\mat{X}$}};
    \tnLeafEE[start at=nodematrix.east, type=e, label={\fontTwo{$I_2$}}]{}{}{}{};
    \tnLeafWW[start at=nodematrix.west, type=w, label={\fontTwo{$I_1$}}]{}{}{}{};
    \node[anchor=north] (dimsmat) at  ($(matrix.south) + (0,-0.5)$) {\fontTwo{$(I_1 \times I_2)$}};

    \Tensor[start at=vector.east,
    xshift=17,
    width=\J,
    height=\I,
    depth=\K,
    anchor=west,
    colorf=frontal,
    colorl=frontal,
    colorh=frontal,
    name=tensor]{};
    \node[anchor=west] (equal4) at ($(tensor.east) + (3,0)$) {\fontOne{$=$}};
    \node[tnode, outer color=frontal, anchor=west] (nodetensor) at ($(equal4.east) + (6.2,0)$){};
    \node[anchor=center] (xtensor) at  ($(tensor.center) + (0,0)$) {\fontOne{$\ten{X}$}};
    \tnLeafNN[start at=nodetensor.north, type=n, label=\fontTwo{$I_1$}]{}{}{}{};
    \tnLeafSE[start at=nodetensor.south east, type=se, label=\fontTwo{$I_3$}]{}{}{}{};
    \tnLeafSW[start at=nodetensor.south west, type=sw, label=\fontTwo{$I_2$}]{}{}{}{};
    \node[anchor=north] (dimsmat) at  ($(tensor.south) + (0,-0.5)$) {\fontTwo{$(I_1 \times I_2 \times I_3)$}};

\end{tikzpicture}

%% file: 04-tikz/basic-operations-tn-rep.tex

\newcommand{\colorten}{frontal} 
\newcommand{\leafsize}{7cm}
\newcommand{\eqshift}{7cm}
\newcommand{\sizeten}{2.5cm}

\newcommand{\yshift}{8.5cm} 
\newcommand{\ylshift}{0.8} 

\newcommand{\fThree}[1]{\scalebox{6.5}{#1}} 
\newcommand{\fontTwo}[1]{\scalebox{8.5}{#1}} 
\newcommand{\fontFour}[1]{\scalebox{5}{#1}} 

\begin{tikzpicture}
\tikzset{
    core/.style={
        circle,
        outer color=\colorten,
        inner color=white,
        no sep,
        minimum size=\sizeten,
        anchor=north,
        draw,
        line width=3pt,
    },
    myleafh/.style={
        rectangle,
        no sep,
        minimum height=0pt,
        minimum width=\leafsize,
        draw,
        line width=3pt,
        anchor=west
    },
    myleafv/.style={
        rectangle,
        no sep,
        minimum height=\leafsize,
        minimum width=0pt,
        draw,
        line width=3pt,
        anchor=north
    },
    myleafvtkd/.style={
        rectangle,
        no sep,
        minimum height=\leafsizetkd,
        minimum width=0pt,
        draw,
        line width=3pt,
        anchor=north
    },
    myline/.style={
        rectangle,
        no sep,
        minimum height=0pt,
        minimum width=30pt,
        draw,
        line width=2pt,
        rotate=45,
        anchor=center
    }
}
\node[core] (A) at (0,0) {};
\node[myleafh, anchor=east] (A_l1) at (A.west) {};
\node[myleafh, anchor=west] (A_l2) at (A.east) {};
\node[myline] () at (A_l1) {};
\node[myline] () at (A_l2) {};
\node[anchor=south] () at ($(A.north) + (0,\ylshift)$) {\fThree{$\mat{X}$}};
\node[anchor=south] () at ($(A_l1.north) + (0,\ylshift)$) {\fontFour{$I_1$}};
\node[anchor=south] () at ($(A_l2.north) + (0,\ylshift)$) {\fontFour{$I_2$}};

\node[core] (A_1) at ($(A.south) + (0,-\yshift)$) {};
\node[myleafh, anchor=east] (A_1_l1) at (A_1.west) {};
\node[myleafh, anchor=west] (A_1_l2) at (A_1.east) {};
\node[myline] () at (A_1_l1) {};
\node[myline] () at (A_1_l2) {};
\node[anchor=south] () at ($(A_1.north) + (0,\ylshift)$) {\fThree{$\mat{X}$}};
\node[anchor=south] () at ($(A_1_l1.north) + (0,\ylshift)$) {\fontFour{$I_1$}};
\node[anchor=south] () at ($(A_1_l2.north) + (0,\ylshift)$) {\fontFour{$I_2=J_1$}};
\node[core, anchor=west] (B_1) at ($(A_1_l2.east)$) {};
\node[myleafh, anchor=west] (B_1_l2) at (B_1.east) {};
\node[myline] () at (B_1_l2) {};
\node[anchor=south] () at ($(B_1.north) + (0,\ylshift)$) {\fThree{$\mat{Y}$}};
\node[anchor=south] () at ($(B_1_l2.north) + (0,\ylshift)$) {\fontFour{$J_2$}};

\node[core] (A_2) at ($(A_1.south) + (0,-\yshift)$) {};
\node[myleafh, anchor=east, rotate=-45] (A_2_l1) at (A_2.north west) {};
\node[myleafh, anchor=west] (A_2_l2) at (A_2.east) {};
\node[myleafh, anchor=east, rotate=45] (A_2_l3) at (A_2.south west) {};
\node[myline] () at (A_2_l1) {};
\node[myline] () at (A_2_l2) {};
\node[myline,rotate=90] () at (A_2_l3) {};
\node[anchor=south] () at ($(A_2.north) + (0,\ylshift)$) {\fThree{$\ten{X}$}};
\node[anchor=south] () at ($(A_2_l1.north) + (0,\ylshift)$) {\fontFour{$I_1$}};
\node[anchor=south] () at ($(A_2_l2.north) + (0,\ylshift)$) {\fontFour{$I_3=J_1$}};
\node[anchor=north] () at ($(A_2_l3.south) + (0,-\ylshift)$) {\fontFour{$I_2$}};
\node[core, anchor=west] (B_2) at ($(A_2_l2.east)$) {};
\node[myleafh, anchor=west, rotate=45] (B_2_l2) at (B_2.north east) {};
\node[myleafh, anchor=west, rotate=-45] (B_2_l3) at (B_2.south east) {};
\node[myline, rotate=90] () at (B_2_l2) {};
\node[myline, rotate=0] () at (B_2_l3) {};
\node[anchor=south] () at ($(B_2.north) + (0,\ylshift)$) {\fThree{$\ten{Y}$}};
\node[anchor=south] () at ($(B_2_l2.north) + (0,\ylshift)$) {\fontFour{$J_3$}};
\node[anchor=north] () at ($(B_2_l3.south) + (0,-\ylshift)$) {\fontFour{$J_2$}};

\node[core] (A_3) at ($(A_2.south) + (0,-\yshift)$) {};
\node[myleafh, anchor=west] (A_3_l2) at (A_3.east) {};
\node[core, anchor=west] (B_3) at ($(A_3_l2.east)$) {};
\draw[bend left, line width=3pt , -]  (A_3.north) to node [midway,anchor=center] (A_3_l1) {} (B_3.north);
\draw[bend right, line width=3pt, -]  (A_3.south) to node [midway,anchor=center] (A_3_l3) {} (B_3.south);
\node[myline] () at (A_3_l1) {};
\node[myline] () at (A_3_l2) {};
\node[myline] () at (A_3_l3) {};
\node[anchor=south] () at ($(A_3_l1.north) + (0,\ylshift)$) {\fontFour{$I_1=J_1$}};
\node[anchor=south] () at ($(A_3_l2.north) + (0,\ylshift)$) {\fontFour{$I_2=J_2$}};
\node[anchor=north] () at ($(A_3_l3.south) + (0,-\ylshift)$) {\fontFour{$I_3=J_3$}};
\node[anchor=east] () at ($(A_3.west) + (-\ylshift,0)$) {\fThree{$\ten{X}$}};
\node[anchor=west] () at ($(B_3.east) + (\ylshift,0)$) {\fThree{$\ten{Y}$}};

\node[core] (X) at ($(A_3.south) + (0,-\yshift)$) {};
\node[anchor=south] () at ($(X.north) + (0,\ylshift)$) {\fThree{$\ten{X}$}};
\node[myleafh, anchor=east] (X_l1) at (X.west){};
\node[myline, rotate=0] () at (X_l1) {};
\node[anchor=south] () at ($(X_l1.north) + (0,\ylshift)$) {\fontFour{$I_1$}};
\ExtractCoordinateX{B_3}
\ExtractCoordinateY{X}
\node[unfolding, minimum size=1.5cm,anchor=center] (kr) at (\XCoord,\YCoord) {};
\node[anchor=east] () at ($(kr.west) + (-\ylshift,0)$) {\fontFour{ $1$}};
\draw[bend left, line width=3pt , -]  (X.north east) to node [midway,anchor=center] (A_3_l1) {} (kr.north);
\draw[bend right, line width=3pt, -]  (X.south east) to node [midway,anchor=center] (A_3_l3) {} (kr.south);
\node[myline] () at (A_3_l1) {};
\node[myline] () at (A_3_l3) {};
\node[anchor=south] () at ($(A_3_l1.north) + (0,\ylshift)$) {\fontFour{$I_3$}};
\node[anchor=north] () at ($(A_3_l3.south) + (0,-\ylshift)$) {\fontFour{$I_2$}};
\node[myleafh, anchor=west](kr_l) at (kr.east){};
\node[myline, rotate=0] () at (kr_l) {};
\node[anchor=south] () at ($(kr_l.north) + (0,\ylshift)$) {\fontFour{$I_2 I_3$}};

\node[anchor=center] (eq1) at ($(A_l2.east) + (\sizeten,0) + (\leafsize,0) + (\eqshift,0)$) {\fontTwo{ $\Leftrightarrow$}};
\node[anchor=center] (eq2) at ($(B_1_l2.east) + (\eqshift,0)$) {\fontTwo{ $=$}};
\node[anchor=center] (eq3) at ($(B_2.east) + (\leafsize,0) + (\eqshift,0)$) {\fontTwo{ $=$}};
\node[anchor=center] (eq4) at ($(B_3.east) + (\leafsize,0) + (\eqshift,0)$) {\fontTwo{ $=$}};
\ExtractCoordinateX{eq4}
\ExtractCoordinateY{kr}
\node[anchor=center] (eq5) at (\XCoord,\YCoord) {\fontTwo{ $\Leftrightarrow$}};

\node[myleafh, anchor=west] (A_T_l1) at ($(eq1.center) + (\eqshift,0)$) {};
\node[core,anchor=west] (A_T) at (A_T_l1.east) {};
\node[myleafh, anchor=west] (A_T_l2) at (A_T.east) {};
\node[myline] () at (A_T_l1) {};
\node[myline] () at (A_T_l2) {};
\node[anchor=south] () at ($(A_T_l1.north) + (0,\ylshift)$) {\fontFour{$I_2$}};
\node[anchor=south] () at ($(A_T_l2.north) + (0,\ylshift)$) {\fontFour{$I_1$}};
\node[anchor=south] () at ($(A_T.north) + (0,\ylshift)$) {\fThree{$\mat{X}^T$}};

\node[myleafh, anchor=west] (C_l1) at ($(eq2.center) + (\eqshift,0)$) {};
\node[core,anchor=west] (C) at (C_l1.east) {};
\node[myleafh, anchor=west] (C_l2) at (C.east) {};
\node[myline] () at (C_l1) {};
\node[myline] () at (C_l2) {};
\node[anchor=south] () at ($(C_l1.north) + (0,\ylshift)$) {\fontFour{$I_1$}};
\node[anchor=south] () at ($(C_l2.north) + (0,\ylshift)$) {\fontFour{$J_2$}};
\node[anchor=south] () at ($(C.north) + (0,\ylshift)$) {\fThree{$\mat{Z} = \mat{XY}$}};

\node[core,anchor=west] (C_2) at ($(eq3.center) + (\eqshift,0) + (\leafsize,0)$) {};
\node[myleafh, anchor=east, rotate=-45] (C_2_l1) at (C_2.north west) {};
\node[myleafh, anchor=west, rotate=45] (C_2_l2) at (C_2.north east) {};
\node[myleafh, anchor=west, rotate=-45] (C_2_l3) at (C_2.south east) {};
\node[myleafh, anchor=east, rotate=45] (C_2_l4) at (C_2.south west) {};
\node[myline] () at (C_2_l1) {};
\node[myline, rotate=90] () at (C_2_l2) {};
\node[myline, rotate=0] () at (C_2_l3) {};
\node[myline,rotate=90] () at (C_2_l4) {};
\node[anchor=south] () at ($(C_2.north) + (0,\ylshift)$) {\fThree{$\ten{Z}$}};
\node[anchor=south] () at ($(C_2_l1.north) + (0,\ylshift)$) {\fontFour{$I_1$}};
\node[anchor=south] () at ($(C_2_l2.north) + (0,\ylshift)$) {\fontFour{$J_3$}};
\node[anchor=north] () at ($(C_2_l3.south) + (0,-\ylshift)$) {\fontFour{$J_2$}};
\node[anchor=north] () at ($(C_2_l4.south) + (0,-\ylshift)$) {\fontFour{$I_2$}};

\node[core,anchor=west] (C_3) at ($(eq4.center) + (\eqshift,0) + (\leafsize,0)$) {};
\node[anchor=south] () at ($(C_3.north) + (0,\ylshift)$) {\fThree{$z$}};

\node[myleafh, anchor=west] (A_T_l1) at ($(eq5.center) + (\eqshift,0)$) {};
\node[core,anchor=west] (A_T) at (A_T_l1.east) {};
\node[myleafh, anchor=west] (A_T_l2) at (A_T.east) {};
\node[myline] () at (A_T_l1) {};
\node[myline] () at (A_T_l2) {};
\node[anchor=south] () at ($(A_T_l1.north) + (0,\ylshift)$) {\fontFour{$I_1$}};
\node[anchor=south] () at ($(A_T_l2.north) + (0,\ylshift)$) {\fontFour{$I_2 I_3$}};
\node[anchor=south] () at ($(A_T.north) + (0,\ylshift)$) {\fThree{$\unfold{X}{1}$}};

\end{tikzpicture}

%% file: 04-tikz/TT.tex

\newcommand{\sizeIone}{7}
\newcommand{\sizeItwo}{7}
\newcommand{\sizeIthree}{7}
\newcommand{\sizeR}{4}

\newcommand{\fontOne}[1]{\scalebox{10.5}{#1}} 
\newcommand{\fontTwo}[1]{\scalebox{5.5}{#1}} 
\newcommand{\fontThree}{\fontsize{50}{70}\selectfont}

\newcommand{\shiftA}{1.5}
\newcommand{\shiftSizes}{0.5}
\newcommand{\eqshift}{3}
\newcommand{\rhsShift}{2}
\newcommand{\coreShift}{\rhsShift + 1}
\newcommand{\lhsShift}{2}
\newcommand{\tShift}{0.125}

\begin{tikzpicture}
    \coordinate (start) at (0,0);
    \node[anchor=center] (eq1) at (start) {\fontOne{$=$}};

    \Tensor[width=\sizeR,
            height=\sizeR,
            depth=\sizeR,
            start at=start.west,
            anchor=east,
            xshift=-\eqshift-2,
            yshift=0,
            zshift=0,
            name=T23]{};

    \Tensor[width=\sizeR,
            height=\sizeR,
            depth=\sizeR,
            start at=T23.north,
            anchor=south,
            xshift=-0,
            yshift=\lhsShift+\tShift,
            zshift=0,
            name=T13]{};

    \Tensor[width=\sizeR,
            height=\sizeR,
            depth=\sizeR,
            start at=T23.south,
            anchor=north,
            xshift=-0,
            yshift=-\lhsShift-\tShift,
            zshift=0,
            name=T33]{};

    \Tensor[width=\sizeR,
            height=\sizeR,
            depth=\sizeR,
            start at=T13.west,
            anchor=east,
            xshift=-\lhsShift-\tShift,
            yshift=0,
            zshift=0,
            name=T12]{};

    \Tensor[width=\sizeR,
            height=\sizeR,
            depth=\sizeR,
            start at=T23.west,
            anchor=east,
            xshift=-\lhsShift-\tShift,
            yshift=0,
            zshift=0,
            name=T22]{};

    \Tensor[width=\sizeR,
            height=\sizeR,
            depth=\sizeR,
            start at=T33.west,
            anchor=east,
            xshift=-\lhsShift-\tShift,
            yshift=0,
            zshift=0,
            name=T32]{};

    \Tensor[width=\sizeR,
            height=\sizeR,
            depth=\sizeR,
            start at=T12.west,
            anchor=east,
            xshift=-\lhsShift-\tShift,
            yshift=0,
            zshift=0,
            name=T11]{};

    \Tensor[width=\sizeR,
            height=\sizeR,
            depth=\sizeR,
            start at=T22.west,
            anchor=east,
            xshift=-\lhsShift-\tShift,
            yshift=0,
            zshift=0,
            name=T21]{};

    \Tensor[width=\sizeR,
            height=\sizeR,
            depth=\sizeR,
            start at=T32.west,
            anchor=east,
            xshift=-\lhsShift-\tShift,
            yshift=0,
            zshift=0,
            name=T31]{};

    \node[anchor=center]
        () at ($(T32.south)+(0,-2)$)
        {\fontTwo{ $\ten{X} \in \R^{I_1 \times I_2 \times I_3 \times I_4 \times I_5}$}};

    \draw
        [decoration={brace, amplitude=20pt, raise=55pt},line width = 3pt, decorate, color=SeaGreen]
        ($(T11.north west) + (1.5, 0)$) --
        ($(T13.north east) + (1.5, 0)$)
        node [above=85pt,midway,sloped] {\fontTwo{\bf{5-th dimension}}};

    \draw
        [decoration={brace, amplitude=20pt, raise=10pt},line width = 3pt, decorate, color=SkyBlue3]
        (T31.south west) --
        (T11.north west)
        node [above=40pt,midway,sloped] {\fontTwo {\bf{4-th dimension}}};


    \fslice[width=\sizeR,
            height =\sizeIone,
            start at=eq1.east,
            anchor=west,
            xshift=\eqshift,
            yshift=0,
            zshift=0,
            color=colorFrontal,
            name=G1]{};
    \node[anchor=center] at (G1.center) {\fontTwo{$\mat{A}$}};
    \node[anchor=north] at ($(G1.south)+(0,-\shiftSizes)$) {\fontThree $(I_1 \times R_1)$};

    \Tensor[width=\sizeR,
            height=\sizeR,
            depth=\sizeR,
            start at=G1.north east,
            anchor=north west,
            xshift=\rhsShift,
            yshift=0,
            zshift=0,
            colorf=colorLateral,
            colorl=colorLateral,
            colorh=colorLateral,
            name=G2]{};
    \node[anchor=center] () at ($(G2.center)+(0,0)$) {\fontTwo{ $\ten{G}^{(2)}$}};
    \node[anchor=north] at ($(G2.south)+(0,-\shiftSizes)$) {\fontThree $(R_1 \times I_2 \times R_2)$};

    \Tensor[width=\sizeR,
            height=\sizeR,
            depth=\sizeR,
            start at=G2.north east,
            anchor=north west,
            xshift=\coreShift,
            yshift=0,
            zshift=0,
            colorf=colorHorizontal,
            colorl=colorHorizontal,
            colorh=colorHorizontal,
            name=G3]{};
    \node[anchor=center] () at ($(G3.center)+(0,0)$) {\fontTwo{ $\ten{G}^{(3)}$}};
    \node[anchor=north] at ($(G3.south)+(0,-\shiftSizes)$) {\fontThree $(R_2 \times I_3 \times R_3)$};

    \Tensor[width=\sizeR,
            height=\sizeR,
            depth=\sizeR,
            start at=G3.north east,
            anchor=north west,
            xshift=\coreShift,
            yshift=0,
            zshift=0,
            colorf=SkyBlue3,
            colorl=SkyBlue3,
            colorh=SkyBlue3,
            name=G4]{};
    \node[anchor=center] () at ($(G4.center)+(0,0)$) {\fontTwo{ $\ten{G}^{(4)}$}};
    \node[anchor=north] at ($(G4.south)+(0,-\shiftSizes)$) {\fontThree $(R_3 \times I_4 \times R_4)$};

    \fslice[width=\sizeIone,
            height =\sizeR,
            start at=G4.north east,
            anchor=north west,
            xshift=\coreShift,
            yshift=0,
            zshift=0,
            color=SeaGreen,
            name=G5]{};
    \node[anchor=center] at (G5.center) {\fontTwo{$\mat{B}$}};
    \node[anchor=north] at ($(G5.south)+(0,-\shiftSizes)$) {\fontThree $(R_4 \times I_5)$};

\end{tikzpicture}

%% file: 04-tikz/TT-rep.tex

\newcommand{\fontOne}[1]{\scalebox{10.5}{#1}} 
\newcommand{\fontTwo}[1]{\scalebox{4}{#1}} 
\newcommand{\fontThree}[1]{\scalebox{3}{#1}} 
\newcommand{\nodeSize}[0]{3cm}
\newcommand{\eqshift}{10cm}

\begin{tikzpicture}
    \setTNsTikzset[leaf size=8cm]

    \coordinate (start) at (0,0);
    \node[anchor=center] (eq1) at (start) {\fontOne{$=$}};

    \tnode[name=X, start at=eq1.west, xshift=-\eqshift]{outer color=horizontal}{\fontTwo{$\ten{X}$}};
    \tnLeafSW[start at=X, label=\fontThree{$I_1$}]{}{}{}{};
    \tnLeafSE[start at=X, label=\fontThree{$I_2$}]{}{}{}{};
    \tnLeafNW[start at=X, label=\fontThree{$I_3$}]{}{}{}{};
    \tnLeafNN[start at=X, label=\fontThree{$I_4$}]{}{}{}{};
    \tnLeafNE[start at=X, label=\fontThree{$I_5$}]{}{}{}{};

    \mnode[name=G1, start at=eq1.west, xshift=\eqshift]{minimum size=\nodeSize}{};
    \node[anchor=south] at (G1.north) {\fontTwo{$\mat{A} = \ten{G}^{(1)}$}};
    \tnLeafSS[start at=G1, label=\fontThree{$I_1$}]{}{}{}{};
    \tnLeafEE[name=r1, start at=G1, label=\fontThree{$R_1$}]{}{}{}{};

    \tnode[name=G2, start at=r1.east]{minimum size=\nodeSize}{};
    \node[anchor=south] at (G2.north) {\fontTwo{$\ten{G}^{(2)}$}};
    \tnLeafSS[start at=G2, label=\fontThree{$I_2$}]{}{}{}{};
    \tnLeafEE[name=r2, start at=G2, label=\fontThree{$R_2$}]{}{}{}{};

    \tnode[name=G3, start at=r2.east]{minimum size=\nodeSize}{};
    \node[anchor=south] at (G3.north) {\fontTwo{$\ten{G}^{(3)}$}};
    \tnLeafSS[start at=G3, label=\fontThree{$I_3$}]{}{}{}{};
    \tnLeafEE[name=r3, start at=G3, label=\fontThree{$R_3$}]{}{}{}{};

    \tnode[name=G4, start at=r3.east]{minimum size=\nodeSize}{};
    \node[anchor=south] at (G4.north) {\fontTwo{$\ten{G}^{(4)}$}};
    \tnLeafSS[start at=G4, label=\fontThree{$I_{4}$}]{}{}{}{};
    \tnLeafEE[name=rn-1, start at=G4, label=\fontThree{$R_{4}$}]{}{}{}{};

    \mnode[name=G5, start at=rn-1.east]{minimum size=\nodeSize}{};
    \node[anchor=south] at (G5.north) {\fontTwo{$\mat{B} = \ten{G}^{(5)}$}};
    \tnLeafSS[start at=G5, label=\fontThree{$I_{5}$}]{}{}{}{};

\end{tikzpicture}

%% file: 04-tikz/TTCP-rep.tex

\newcommand{\localLeafSize}{8.5cm} 
\newcommand{\localOriginalTensorMult}{0.8} 
\newcommand{\localNodeSizes}{2cm} 
\newcommand{\xShift}{10cm}  
\newcommand{\yShiftTT}{2cm} 

\newcommand{\fontLeafs}[1]{\scalebox{4}{#1}} 
\newcommand{\fontMNodes}[1]{\scalebox{5}{#1}} 
\newcommand{\fontEq}[1]{\scalebox{10}{#1}} 
\newcommand{\fontText}[1]{\scalebox{5}{#1}} 

\newcommand{\focusColor}[0]{horizontal}
\newcommand{\focusLineSize}[0]{4pt}

\begin{tikzpicture}
\tikzset{
    inline text/.style={
        anchor=north,
        inner sep=15pt
    },
    inner spacing/.style={
        inner sep=15pt
    },
    mydash/.style={
        line width=8pt,
        dashed,
        color=\focusColor,
        dash pattern=on 15pt off 15pt
    },
}
\setTNsTikzset[%
    leaf size=\localLeafSize,%
    mnode size=\localNodeSizes,%
    tnode size=\localNodeSizes,%
]
\coordinate (start) at (0,0);

\mnode[name=node, start at=start, xshift=-\xShift]{}{};
\node[anchor=south] at (node.north) {\fontMNodes{$\ten{G}_y^{(3)}=\mat{B}_y$}};
\tnLeafSS[name=leaf-ss, start at=node, label=\fontLeafs{$J_3$}]{}{}{}{};
\tnLeafWW[name=leaf-ww, start at=node, label=\fontLeafs{$P_2$}]{}{}{anchor=north}{};
\node (ref0) at (leaf-ss.south) {};

\tnode[name=node, start at=leaf-ww.west]{}{};
\node[anchor=south] at (node.north) {\fontMNodes{$\ten{G}_y^{(2)}$}};
\tnLeafSS[start at=node, label=\fontLeafs{$J_2 = I_3$}]{minimum width=\focusLineSize, fill=\focusColor, color=\focusColor}{color=\focusColor}{}{};
\tnLeafWW[name=leaf-ww, start at=node, label=\fontLeafs{$P_1$}]{}{}{anchor=north}{};

\mnode[name=node, start at=leaf-ww.west]{}{};
\node[anchor=south] at (node.north) {\fontMNodes{$\ten{G}_y^{(1)} = \mat{A}_y$}};
\tnLeafSS[start at=node, label=\fontLeafs{$J_1$}]{}{}{}{};

\mnode[name=node, start at=ref0, yshift=-\localLeafSize-\yShiftTT]{}{};
\node[anchor=north] at (node.south) {\fontMNodes{$\ten{G}_x^{(4)} = \mat{B}_x$}};
\tnLeafNN[name=leaf-ss, start at=node, label=\fontLeafs{$I_4$}]{}{}{}{};
\tnLeafWW[name=leaf-ww, start at=node, label=\fontLeafs{$R_3$}]{}{}{}{};

\tnode[name=node, start at=leaf-ww.west]{}{};
\node[anchor=north] at (node.south) {\fontMNodes{$\ten{G}_x^{(3)}$}};
\tnLeafNN[start at=node, label=\fontLeafs{$I_3 = J_2$}]{minimum width=\focusLineSize, fill=\focusColor, color=\focusColor}{color=\focusColor}{}{};
\tnLeafWW[name=leaf-ww, start at=node, label=\fontLeafs{$R_2$}]{}{}{}{};

\tnode[name=node, start at=leaf-ww.west]{}{};
\node[anchor=north] at (node.south) {\fontMNodes{$\ten{G}_x^{(2)}$}};
\tnLeafNN[start at=node, label=\fontLeafs{$I_2$}]{}{}{}{};
\tnLeafWW[name=leaf-ww, start at=node, label=\fontLeafs{$R_1$}]{}{}{}{};

\mnode[name=node, start at=leaf-ww.west]{}{};
\node[anchor=north] at (node.south) {\fontMNodes{$\ten{G}_x^{(1)} = \mat{A}_x$}};
\tnLeafNN[start at=node, label=\fontLeafs{$I_1$}]{}{}{}{};

\tnode[name=node, start at=node, xshift=-1.5*\localLeafSize]{}{};
\node[anchor=south] at (node.north) {\fontMNodes{$\ten{X}$}};
\tnLeafNE[start at=node, label=\fontLeafs{$I_1$}]{}{}{}{\localOriginalTensorMult*\localLeafSize};
\tnLeafSW[start at=node, label=\fontLeafs{$I_2$}]{}{}{}{\localOriginalTensorMult*\localLeafSize};
\tnLeafSE[start at=node, label=\fontLeafs{$I_3 = J_2$}]{minimum height=\focusLineSize, fill=\focusColor, color=\focusColor}{color=\focusColor}{}{\localOriginalTensorMult*\localLeafSize};
\tnLeafNW[start at=node, label=\fontLeafs{$I_4$}]{}{}{}{\localOriginalTensorMult*\localLeafSize};
\node (ref5) at (node) {};

\tnode[name=node, start at=node, yshift=2*\localLeafSize+\yShiftTT]{}{};
\node[anchor=south west] at (node.north east) {\fontMNodes{$\ten{Y}$}};
\tnLeafNN[start at=node, label=\fontLeafs{$J_1$}]{}{}{}{\localOriginalTensorMult*\localLeafSize};
\tnLeafSW[start at=node, label=\fontLeafs{$J_2 = I_3$}]{minimum height=\focusLineSize, fill=\focusColor, color=\focusColor}{color=\focusColor}{}{\localOriginalTensorMult*\localLeafSize};
\tnLeafSE[start at=node, label=\fontLeafs{$J_3$}]{}{}{}{\localOriginalTensorMult*\localLeafSize};

\mnode[name=node, start at=start, xshift=\xShift]{}{};
\node[anchor=south] at (node.north) {\fontMNodes{$\ten{G}_y^{(1)}=\mat{A}_y$}};
\tnLeafSS[name=leaf-ss, start at=node, label=\fontLeafs{$J_2 = I_3$}]{minimum width=\focusLineSize, fill=\focusColor, color=\focusColor}{color=\focusColor}{}{};
\tnLeafEE[name=leaf-ee, start at=node, label=\fontLeafs{$P_1$}]{}{}{anchor=north}{};
\node (ref1) at (leaf-ss.south) {};

\tnode[name=node, start at=leaf-ee.east]{}{};
\node[anchor=south] at (node.north) {\fontMNodes{$\ten{G}_y^{(2)}$}};
\tnLeafSS[start at=node, label=\fontLeafs{$J_1$}]{}{}{}{};
\tnLeafEE[name=leaf-ee, start at=node, label=\fontLeafs{$P_2$}]{}{}{anchor=north}{};

\mnode[name=node, start at=leaf-ee.east]{}{};
\node[anchor=south] at (node.north) {\fontMNodes{$\ten{G}_y^{(3)} = \mat{B}_y$}};
\tnLeafSS[start at=node, label=\fontLeafs{$J_3$}]{}{}{}{};

\mnode[name=node, start at=ref1, yshift=-\localLeafSize-\yShiftTT]{}{};
\node[anchor=north] at (node.south) {\fontMNodes{$\ten{G}_x^{(1)} = \mat{A}_x$}};
\tnLeafNN[name=leaf-ss, start at=node, label=\fontLeafs{$I_3 = J_2$}]{minimum width=\focusLineSize, fill=\focusColor, color=\focusColor}{color=\focusColor}{}{};
\tnLeafEE[name=leaf-ee, start at=node, label=\fontLeafs{$R_1$}]{}{}{}{};

\tnode[name=node, start at=leaf-ee.east]{}{};
\node[anchor=north] at (node.south) {\fontMNodes{$\ten{G}_x^{(2)}$}};
\tnLeafNN[start at=node, label=\fontLeafs{$I_2$}]{}{}{}{};
\tnLeafEE[name=leaf-ee, start at=node, label=\fontLeafs{$R_2$}]{}{}{}{};

\tnode[name=node, start at=leaf-ee.east]{}{};
\node[anchor=north] at (node.south) {\fontMNodes{$\ten{G}_x^{(3)}$}};
\tnLeafNN[start at=node, label=\fontLeafs{$I_1$}]{}{}{}{};
\tnLeafEE[name=leaf-ee, start at=node, label=\fontLeafs{$R_3$}]{}{}{}{};

\mnode[name=node, start at=leaf-ee.east]{}{};
\node[anchor=north] at (node.south) {\fontMNodes{$\ten{G}_x^{(4)} = \mat{B}_x$}};
\tnLeafNN[start at=node, label=\fontLeafs{$I_4$}]{}{}{}{};

\tnode[name=node, start at=node, xshift=1.5*\localLeafSize]{}{};
\node[anchor=south] at (node.north) {\fontMNodes{$\ten{X}^\rho$}};
\tnLeafNE[start at=node, label=\fontLeafs{$I_3 = J_2$}]{minimum width=\focusLineSize, fill=\focusColor, color=\focusColor}{color=\focusColor}{}{\localOriginalTensorMult*\localLeafSize};
\tnLeafSE[start at=node, label=\fontLeafs{$I_1$}]{}{}{}{\localOriginalTensorMult*\localLeafSize};
\tnLeafSW[start at=node, label=\fontLeafs{$I_2$}]{}{}{}{\localOriginalTensorMult*\localLeafSize};
\tnLeafNW[start at=node, label=\fontLeafs{$I_3$}]{}{}{}{\localOriginalTensorMult*\localLeafSize};
\node (ref6) at (node) {};

\tnode[name=node, start at=node, yshift=2*\localLeafSize+\yShiftTT]{}{};
\node[anchor=south west] at (node.north east) {\fontMNodes{$\ten{Y}^\rho$}};
\tnLeafNN[start at=node, label=\fontLeafs{$J_2 = I_3$}]{minimum width=\focusLineSize, fill=\focusColor, color=\focusColor}{color=\focusColor}{}{\localOriginalTensorMult*\localLeafSize};
\tnLeafSE[start at=node, label=\fontLeafs{$J_1$}]{}{}{}{\localOriginalTensorMult*\localLeafSize};
\tnLeafSW[start at=node, label=\fontLeafs{$J_3$}]{}{}{}{\localOriginalTensorMult*\localLeafSize};

\node[anchor=center] (eq1) at ($(ref0)!0.5!(ref1)$) {\fontEq{$\Leftrightarrow$}};
\node[anchor=north] () at (eq1.south)  {\fontText{Permutation}};

\node[anchor=center, rotate=90] (eq2) at ($(eq1) + (0,-1.8*\localLeafSize)$) {\fontEq{$\Leftarrow$}};
\node[inline text] () at (eq2.west)  {\fontText{Contraction operation}};

\mnode[name=node, start at=eq2, yshift=-1.2*\localLeafSize]{}{};
\node (ref2) at (node) {};
\node (ref3) at (node) {};
\node[anchor=south] at (node.north) {\fontMNodes{$\mat{A}_x$}};
\tnLeafEE[name=leaf-ee, start at=node, label=\fontLeafs{$I_3 = J_2$}]{minimum height=\focusLineSize, fill=\focusColor, color=\focusColor}{color=\focusColor}{}{};
\tnLeafWW[name=leaf-ww, start at=node, label=\fontLeafs{$R_1$}]{}{}{}{};

\mnode[name=node, start at=leaf-ee.east]{}{};
\node[anchor=south] at (node.north) {\fontMNodes{$\mat{A}_y$}};
\tnLeafEE[name=leaf-ee, start at=node, label=\fontLeafs{$P_1$}]{}{}{}{};
\node (ref4) at (node) {};

\tnode[name=node, start at=leaf-ee.east]{}{};
\node[anchor=south] at (node.north) {\fontMNodes{$\ten{G}_y^{(2)}$}};
\tnLeafSS[start at=node, label=\fontLeafs{$J_1$}]{}{}{}{};
\tnLeafEE[name=leaf-ee, start at=node, label=\fontLeafs{$R_3$}]{}{}{}{};

\mnode[name=node, start at=leaf-ee.east]{}{};
\node[anchor=south] at (node.north) {\fontMNodes{$\mat{B}_y$}};
\tnLeafSS[start at=node, label=\fontLeafs{$J_3$}]{}{}{}{};

\tnode[name=node, start at=leaf-ww.west]{}{};
\node[anchor=south] at (node.north) {\fontMNodes{$\ten{G}_x^{(2)}$}};
\tnLeafSS[, start at=node, label=\fontLeafs{$I_1$}]{}{}{}{};
\tnLeafWW[name=leaf-ww, start at=node, label=\fontLeafs{$R_2$}]{}{}{}{};

\tnode[name=node, start at=leaf-ww.west]{}{};
\node[anchor=south] at (node.north) {\fontMNodes{$\ten{G}_x^{(3)}$}};
\tnLeafSS[, start at=node, label=\fontLeafs{$I_2$}]{}{}{}{};
\tnLeafWW[name=leaf-ww, start at=node, label=\fontLeafs{$R_3$}]{}{}{}{};

\mnode[name=node, start at=leaf-ww.west]{}{};
\node[anchor=south] at (node.north) {\fontMNodes{$\mat{B}_x$}};
\tnLeafSS[, start at=node, label=\fontLeafs{$I_4$}]{}{}{}{};

\draw[mydash]
    let
        \p1=(ref3),
        \p2=(ref4),
        \n1={atan2(\y2-\y1,\x2-\x1)},
        \n2={veclen(\y2-\y1,\x2-\x1)}
    in
        ($ (ref3)!0.5!(ref4) $)
    ellipse [y radius=\n2/2+3cm, x radius=5cm,rotate=90-\n1];

\node[anchor=center, rotate=90] (eq3) at ($(ref2) + (0,-1*\localLeafSize)$) {\fontEq{$\Leftarrow$}};
\node[inline text] () at (eq3.west)  {\fontText{Matrix multiplication}};

\mnode[name=node, start at=eq3, yshift=-1.2*\localLeafSize, xshift=0.5*\localLeafSize]{outer color=\focusColor}{};
\node (ref2) at (node) {};
\node[anchor=south] at (node.north) {\fontMNodes{$\mat{K}$}};
\tnLeafEE[name=leaf-ee, start at=node, label=\fontLeafs{$P_1$}]{}{}{}{};
\tnLeafWW[name=leaf-ww, start at=node, label=\fontLeafs{$R_1$}]{}{}{}{};

\tnode[name=node, start at=leaf-ee.east]{}{};
\node[anchor=south] at (node.north) {\fontMNodes{$\ten{G}_y^{(2)}$}};
\tnLeafSS[start at=node, label=\fontLeafs{$J_1$}]{}{}{}{};
\tnLeafEE[name=leaf-ee, start at=node, label=\fontLeafs{$P_2$}]{}{}{}{};

\mnode[name=node, start at=leaf-ee.east]{}{};
\node[anchor=south] at (node.north) {\fontMNodes{$\mat{B}_y$}};
\tnLeafSS[start at=node, label=\fontLeafs{$J_3$}]{}{}{}{};
\node (ref8) at (node) {};

\tnode[name=node, start at=leaf-ww.west]{}{};
\node[anchor=south] at (node.north) {\fontMNodes{$\ten{G}_x^{(2)}$}};
\tnLeafSS[, start at=node, label=\fontLeafs{$I_1$}]{}{}{}{};
\tnLeafWW[name=leaf-ww, start at=node, label=\fontLeafs{$R_2$}]{}{}{}{};

\tnode[name=node, start at=leaf-ww.west]{}{};
\node[anchor=south] at (node.north) {\fontMNodes{$\ten{G}_x^{(3)}$}};
\tnLeafSS[, start at=node, label=\fontLeafs{$I_2$}]{}{}{}{};
\tnLeafWW[name=leaf-ww, start at=node, label=\fontLeafs{$R_3$}]{}{}{}{};

\mnode[name=node, start at=leaf-ww.west]{}{};
\node[anchor=south] at (node.north) {\fontMNodes{$\mat{B}_x$}};
\tnLeafSS[, start at=node, label=\fontLeafs{$I_4$}]{}{}{}{};
\node (ref7) at (node) {};

\ExtractCoordinateY{ref7}
\ExtractCoordinateX{ref5}
\node[anchor=center] (t) at ($(\XCoord,\YCoord)$) {};
\tnode[name=node, start at=t]{}{};
\node[anchor=north, inner spacing] at (node.south) {\fontMNodes{$\ten{Z}$}};
\tnLeafNN[start at=node, label=\fontLeafs{$I_1$}]{}{}{}{\localOriginalTensorMult*\localLeafSize};
\tnLeafNE[start at=node, label=\fontLeafs{$I_2$}]{}{}{}{\localOriginalTensorMult*\localLeafSize};
\tnLeafSE[start at=node, label=\fontLeafs{$I_4$}]{}{}{}{\localOriginalTensorMult*\localLeafSize};
\tnLeafSW[start at=node, label=\fontLeafs{$J_1$}]{}{}{}{\localOriginalTensorMult*\localLeafSize};
\tnLeafNW[start at=node, label=\fontLeafs{$J_3$}]{}{}{}{\localOriginalTensorMult*\localLeafSize};
\node[anchor=center] (eq4) at ($(ref7)!0.5!(t)$) {\fontEq{$\Leftarrow$}};
\node[inline text, align=center] () at (eq4.south) {\fontText{Reconstruction}\\[2em] \fontText{with} \\[2em] \fontText{Permutation}};

\ExtractCoordinateY{ref8}
\ExtractCoordinateX{ref6}
\node[anchor=center] (t) at ($(\XCoord,\YCoord)$) {};
\tnode[name=node, start at=t]{}{};
\node[anchor=north, inner spacing] at (node.south) {\fontMNodes{$\ten{Z}$}};
\tnLeafNN[start at=node, label=\fontLeafs{$I_4$}]{}{}{}{\localOriginalTensorMult*\localLeafSize};
\tnLeafNE[start at=node, label=\fontLeafs{$I_2$}]{}{}{}{\localOriginalTensorMult*\localLeafSize};
\tnLeafSE[start at=node, label=\fontLeafs{$I_1$}]{}{}{}{\localOriginalTensorMult*\localLeafSize};
\tnLeafSW[start at=node, label=\fontLeafs{$J_1$}]{}{}{}{\localOriginalTensorMult*\localLeafSize};
\tnLeafNW[start at=node, label=\fontLeafs{$J_3$}]{}{}{}{\localOriginalTensorMult*\localLeafSize};
\node[anchor=center] (eq4) at ($(ref8)!0.5!(t)$) {\fontEq{$\Rightarrow$}};
\node[inline text] () at (eq4.south)  {\fontText{Reconstruction}};

\end{tikzpicture}